\pdfoutput=1
\documentclass[pdflatex,preprint,number]{elsarticle}

\usepackage{amssymb}
\usepackage{amsfonts}
\usepackage{hyperref}
\usepackage{epsfig,amsmath}
\usepackage{tikz,ifthen,calc}
\usepackage{graphicx}
\usepackage{caption}
\usepackage{subcaption}
\usepackage{wrapfig}

\usepackage[margin=1.00 in]{geometry}

\newtheorem{thm}{Theorem}
\newtheorem{lem}[thm]{Lemma}

\newtheorem{propo}[thm]{Proposition}

\newdefinition{rmk}{Remark}
\newdefinition{definition}{Definition}
\newdefinition{example}{Example}
\newproof{pf}{Proof}

\begin{document}

\title{Signed tilings by ribbon $L$ $n$-ominoes, $n$ even, via Gr\"obner bases}

\author[rvt1]{Kenneth Gill}
\ead{KG818321@wcupa.edu}

\author[rvt2]{Viorel Nitica}
\ead{vnitica@wcupa.edu}

\address[rvt1]{Department of Mathematics, West Chester University,
PA 19383, USA}

\address[rvt2]{Department of Mathematics, West Chester University,
PA 19383, USA}

\begin{abstract} Let $\mathcal{T}_n$ be the set of ribbon $L$-shaped $n$-ominoes for some $n\ge 4$ even, and let $\mathcal{T}_n^+$ be $\mathcal{T}_n$ with an extra $2\times 2$ square. We investigate signed tilings of rectangles by $\mathcal{T}_n$ and $\mathcal{T}_n^+$. We show that a rectangle has a signed tiling by $\mathcal{T}_n$ if and only if both sides of the rectangle are even and one of them is divisible by $n$, or if one of the sides is odd and the other side is divisible by $n\left (\frac{n}{2}-2\right ).$ We also show that a rectangle has a signed tiling by $\mathcal{T}_n^+, n\ge 6$ even, if and only if both sides of the rectangle are even, or if one of the sides is odd and the other side is divisible by $n\left (\frac{n}{2}-2\right ).$ Our proofs are based on the exhibition of explicit Gr\"obner bases for the ideals generated by polynomials associated to the tiling sets. In particular, we show that some of the regular tiling results in \emph{ V.~Nitica,  Every tiling of the first quadrant by ribbon $L$ $n$-ominoes follows the rectangular pattern. Open Journal of Discrete Mathematics, {\em 5}, (2015) 11--25,} cannot be obtained from coloring invariants.
\end{abstract}


\date{}
\maketitle

\section{Introduction}

In this article we study tiling problems for regions in a square lattice by certain symmetries of an $L$-shaped polyomino.
Polyominoes were introduced by Golomb in \cite{Golomb1} and the standard reference about this subject is the book \emph{Polyominoes} \cite{Golumb3}. The $L$-shaped polyomino we study is placed in a square lattice and is made out of $n, n\ge 4,$ unit squares, or \emph{cells}. See Figure~\ref{fig:LTetromino1}. In an $a\times b$ rectangle, $a$ is the height and $b$ is the base. We consider translations (only!) of the tiles shown in Figure~\ref{fig:LTetrominoes1}. They are ribbon $L$-shaped $n$-ominoes. A {\em ribbon polyomino}~\cite{Pak} is a simply connected polyomino with no two unit squares lying along a line parallel to the first bisector $y=x$. We denote the set of tiles by $\mathcal{T}_n.$

\begin{figure}[h]
~~~~~~~~~
\begin{subfigure}{.15\textwidth}
\begin{tikzpicture}[scale=.4]
\draw [line width = 1](2.5,1)--(1,1)--(1,2)--(0,2)--(0,0)--(2.5,0);
\draw[dotted, line width=1] (2.6,.5)--(3.9,.5);
\draw [line width = 1] (4,0)--(6,0)--(6,1)--(4,1);
\draw [line width = 1] (1,0)--(1,1);
\draw [line width = 1] (0,1)--(1,1);
\draw [line width = 1] (2,0)--(2,1);
\draw [line width = 1] (5,0)--(5,1);
\end{tikzpicture}
\caption{An $L$ $n$-omino with $n$ cells.}
\label{fig:LTetromino1}
\end{subfigure}
~~~~~~~~~
\begin{subfigure}{.5\textwidth}
\begin{tikzpicture}[scale=.4]
\draw [line width = 1](2.5,1)--(1,1)--(1,2)--(0,2)--(0,0)--(2.5,0);
\draw[dotted, line width=1] (2.6,.5)--(3.9,.5);
\draw [line width = 1] (4,0)--(6,0)--(6,1)--(4,1);
\draw [line width = 1] (1,0)--(1,1);
\draw [line width = 1] (0,1)--(1,1);
\draw [line width = 1] (2,0)--(2,1);
\draw [line width = 1] (5,0)--(5,1);
\draw [line width = 1] (10,1)--(8,1)--(8,2)--(10,2);
\draw [dotted, line width=1] (10.1,1.5)--(11.4,1.5);
\draw [line width = 1] (9,1)--(9,2);
\draw [line width = 1] (11.5,1)--(13,1)--(13,0)--(14,0)--(14,2)--(11.5,2);
\draw [line width = 1] (12,1)--(12,2);
\draw [line width = 1] (13,1)--(13,2);
\draw [line width = 1] (13,1)--(14,1);

\draw [line width = 1] (-9,2)--(-9,4)--(-8,4)--(-8,2);
\draw [dotted, line width=1] (-8.5,1.9)--(-8.5,.6);
\draw [line width = 1] (-9,.5)--(-9,-2)--(-7,-2)--(-7,-1)--(-8,-1)--(-8,.5);
\draw [line width = 1] (-9,3)--(-8,3);
\draw [line width = 1] (-9,0)--(-8,0);
\draw [line width = 1] (-9,-1)--(-8,-1);
\draw [line width = 1] (-8,-1)--(-8,-2);

\draw [line width = 1] (-3,1.5)--(-3,3)--(-4,3)--(-4,4)--(-2,4)--(-2,1.5);
\draw [line width = 1] (-3,2)--(-2,2);
\draw [line width = 1] (-3,3)--(-2,3);
\draw [line width = 1] (-3,3)--(-3,4);
\draw [dotted, line width=1] (-2.5,1.4)--(-2.5,.1);
\draw [line width = 1] (-3,0)--(-3,-2)--(-2,-2)--(-2,0);
\draw [line width = 1] (-3,-1)--(-2,-1);
\end{tikzpicture}
\caption{The set of tiles $\mathcal{T}_n$.}
\label{fig:LTetrominoes1}
\end{subfigure}
\caption{}
\end{figure}

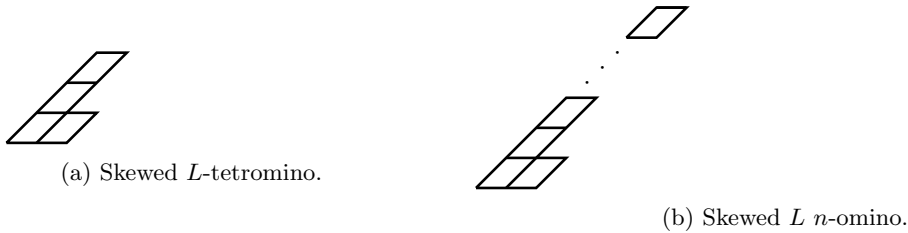
\begin{figure}[h]
~~~~~~~~~
\begin{subfigure}{.30\textwidth}
\begin{tikzpicture}[scale=.4]
\draw [line width = 1](0,0)--(3,3)--(4,3)--(2,1)--(3,1)--(2,0)--(0,0);
\draw [line width = 1](1,1)--(2,1);
\draw [line width = 1](2,2)--(3,2);
\draw [line width = 1](1,0)--(2,1);
\end{tikzpicture}
\caption{Skewed $L$-tetromino.}
\label{fig:skewed-LTetromino1}
\end{subfigure}
~~~~~~~~~
\begin{subfigure}{.5\textwidth}
\begin{tikzpicture}[scale=.4]
\draw [line width = 1](0,0)--(3,3)--(4,3)--(2,1)--(3,1)--(2,0)--(0,0);
\draw [line width = 1](1,1)--(2,1);
\draw [line width = 1](2,2)--(3,2);
\draw [line width = 1](1,0)--(2,1);

\draw [line width = 1] (5,5)--(6,6)--(7,6)--(6,5)--(5,5);
\draw [line width = 1] (6,6)--(7,6);
\node at (3.7, 3.5) {$\cdot$};
\node at (4.2, 4) {$\cdot$};
\node at (4.7, 4.5) {$\cdot$};
\end{tikzpicture}
\caption{Skewed $L$ $n$-omino.}
\label{fig:skewd-LTetrominoesn}
\end{subfigure}
\caption{Skewed polyominoes}
\end{figure}

Related papers are~\cite{CLNS},~\cite{nitica-L-shaped}, investigating tilings by $\mathcal{T}_n, n$ even. In \cite{CLNS} we look at tilings by $\mathcal{T}_n$ in the particular case $n=4$. The starting point was a problem from recreational mathematics. We recall that a replicating tile is one that can make larger copies of itself. The order of replication is the number of initial tiles that fit in the larger copy. Replicating tiles were introduced by Golomb in~\cite{Golomb-rep}. In~\cite{Nitica} we study replication of higher orders for several tiles introduced in~\cite{Golomb-rep}. In particular, we suggested that the skewed $L$-tetromino showed in Figure~\ref{fig:skewed-LTetromino1} is not replicating of order $k^2$ for any odd $k$. The question is equivalent to that of tiling a $k$-inflated copy of the straight $L$-tetromino using only the ribbon orientations of an $L$-tetromino. The question is solved in~\cite{CLNS}, where it is shown that the $L$-tetromino is not replicating of any odd order. This is a consequence of a stronger result: a tiling of the first quadrant by $\mathcal{T}_4$ always follows the rectangular pattern, that is, the tiling reduces to a tiling by $4\times 2$ and $2\times 4$ rectangles, each tiled in turn by two tiles from $\mathcal{T}_4$.

The results in~\cite{CLNS} are generalized in~\cite{nitica-L-shaped} to $\mathcal{T}_n, n$ even. The main result shows that any tiling of the first quadrant by $\mathcal{T}_n$ reduces to a tiling by $2\times n$ and $n\times 2$ rectangles. An application is the characterization of all rectangles that can be tiled by $\mathcal{T}_n, n$ even: a rectangle can be tiled by $\mathcal{T}_n, n$ even, if and only if both sides are even and at least one side is divisible by $n$. The rectangular pattern persists if one adds an extra $2\times 2$ tile to $\mathcal{T}_n, n$ even. The new tiling set is denoted by $\mathcal{T}_n^+$. A rectangle can be tiled by $\mathcal{T}_n^+$ if and only if it has both sides even. The main result also implies that a skewed $L$-shaped $n$-omino, $n$ even, (see Figure~\ref{fig:skewd-LTetrominoesn}) is not a replicating tile of order $k^2$ for any odd $k$.
This development shows that the limitation of the orientations of the tiles can be of interest, in particular when investigating tiling problems in a skewed lattice.

Signed tilings (see~\cite{con-lag}) are also of interest. These are finite
placements of tiles on a plane, with weights +1 or -1 assigned to each of the tiles. We say that they tile a region R
if the sum of the weights of the tiles is 1 for every cell inside R and 0 for every cell elsewhere. The existence of a regular tiling clearly
implies the existence of a signed tiling. Many times solving a tiling problem can be reduced to a coloring argument. It was shown in~\cite{con-lag} that the most general argument of this type is equivalent to the existence of a signed tiling. Consequently, different conditions for regular versus signed tilings can be used to show that certain tiling arguments are stronger then coloring arguments. By looking at signed tilings of rectangles by $\mathcal{T}_n$ and $\mathcal{T}_n^+$, $n$ even, we show that some of the results in~\cite{nitica-L-shaped} cannot be obtained via coloring arguments.

A useful tool in the study of signed tilings is a Gr\"obner basis associated to the polynomial ideal generated by the tiling set. See Bodini and Nouvel~\cite{B-N}. One can associate to any cell in the square lattice a monomial in two variable $x,y$. If the coordinates of the lower left corner of the cell are $(\alpha,\beta)$, one associates $x^{\alpha}y^{\beta}$. This correspondence associates to any bounded tile a Laurent polynomial with all coefficients 1. The polynomial associated to a tile $P$ is denoted by $f_P$. The polynomial associated to a tile translated by an integer vector $(\gamma, \delta)$ is the initial polynomial multiplied by the monomial $x^{\gamma}y^{\delta}$. If the region we tile is bounded and the tile set consists of bounded tiles, then the problem can be translated in the first quadrant via a translation by an integer vector, and one can work only with regular polynomials in $\mathbb{Z}[X,Y]$. See Theorem~\ref{thm-groebner-tiling} below.

Signed tilings by ribbon $L$ $n$-ominoes, $n$ odd are studied in~\cite{Nit-odd}, where we show that a rectangle can be signed tiled by ribbon $L$ $n$-ominoes, $n$ odd, if and only if it has a side divisible by $n$.

The main results of the paper are the following:

\begin{thm}\label{thm:main} A rectangle can be signed tiled by $\mathcal{T}_n$, $n\ge 6$ even, if and only if both sides of the rectangle are even and one of them is divisible by $n$, or one of the sides is odd and the other is divisible by $n\left (\frac{n}{2}-2\right ).$
\end{thm}

Theorem~\ref{thm:main} is proved in Section~\ref{s:4}, after finding a Gr\"obner basis for $\mathcal{T}_n$ in Section~\ref{s:3}.
Theorem~\ref{thm:main} shows that some tiling results for $\mathcal{T}_n, n\ge 6$ even, in~\cite{Nit-odd} cannot be found via coloring arguments.
We recall that it is shown in~\cite{CLNS} that a rectangle is signed tiled by $\mathcal{T}_4$ if and only if the sides are even and one side is divisible by 4.

\begin{thm}\label{thm:main++}  A rectangle can be signed tiled by $\mathcal{T}^+_4$ if and only if both sides are even. A rectangle can be signed tiled by $\mathcal{T}^+_n$, $n\ge 6$ even, if and only if it has both sides even or one side is odd and the other side is divisible by $n(\frac{n}{2}-2)$.
\end{thm}

Theorem~\ref{thm:main++} is proved in Section~\ref{s:4-bis}, after finding a Gr\"obner basis for $\mathcal{T}_n^+$ in Section~\ref{s:3-bis}. Theorem~\ref{thm:main++} shows that some tiling results for $\mathcal{T}_n^+, n\ge 6$ even, in~\cite{Nit-odd} cannot be found via coloring arguments.

Due to the Gr\"obner basis that we exhibit for $\mathcal{T}_n, n\ge 6$ even, we also have:

\begin{propo}\label{thm:main-coro93} A $k$-inflated copy of the ribbon $L$ $n$-omino, $n\ge 6$ even, has a signed tiling by $\mathcal{T}_n$ if and only if $k$ is even or $k$ is odd and divisible by $\left (\frac{n}{2}-2\right )$.
\end{propo}

The proof of Proposition~\ref{thm:main-coro93} is shown in Section~\ref{s:5}.

Barnes~\cite{Barnes1, Barnes2} developed a method for solving signed tiling problems with complex number weights, which applied to our tiling sets gives:

\begin{thm}\label{thm-main-barnes} If complex number weights are used, a rectangle can be signed tiled by $\mathcal{T}_n, n\ge 6$ even, if and only if it has a side divisible by $n$. If only integer weights are used, a rectangle that has a side divisible by $n$ and all cells labeled by the same multiple of $\left (\frac{n}{2}-2\right )$ can be signed tiled by $\mathcal{T}_n, n\ge 6$ even.
\end{thm}

Theorem~\ref{thm-main-barnes} is proved in Section~\ref{s:barnes}. A Gr\"obner basis for the tiling set helps even if Barnes method is used.

\begin{thm}\label{t:barnes+} If complex number weights are used, a rectangle can be signed tiled by $\mathcal{T}^+_n, n\ge 6$ even, if and only if it has an even side, and a rectangle can be signed tiled by $\mathcal{T}^+_4$ if and only if both sides are even.
\end{thm}

Theorem~\ref{t:barnes+} is proved in Section~\ref{s:barnes+}. It is not clear to us if last statement in Theorem~\ref{thm-main-barnes} implies Theorem~\ref{thm:main} and if Theorem~\ref{t:barnes+} implies Theorem~\ref{thm:main++}. Guided by the work here and in~\cite{nitica-L-shaped}, we conclude that Gr\"obner basis method for solving signed tiling problems with integer weights is sometimes more versatile and leads to stronger results then Barnes method.

The methods we use in this paper are well known when applied to a particular tiling problem. Here we apply them uniformly to solve an infinite collection of problems. Our hope was to see some regularity in the Gr\"obner bases associated to other infinite families of tiling sets, such as the family $\mathcal{T}_{m,n,p}$ investigated in~\cite{nitica-L-shaped}. We recall that if $m,p$ are odd and $n$ is even, tilings of the first quadrant by this family follow the rectangular pattern.  Nevertheless, our hopes were not validated. The subfamily $\mathcal{T}_{3,n,3}, n$ even, has a wide variety of Gr\"obner bases. Thus for this family we understand regular tilings of rectangles, but cannot decide if the results follow from coloring invariants.

\section{Summary of Gr\"obner basis theory}\label{s:2}

Let $R[\underline{X}] = R[X_1, \dots , X_k]$ be the ring of polynomials with coefficients in a principal ideal domain (PID) $R$. A \emph{term} in the variables $x_1,\dots,x_k$ is a power product $x_1^{\alpha_1}x_2^{\alpha_2}\dots x_{\ell}^{\alpha_{\ell}}$ with $\alpha_i\in \mathbb{N}, 1\le i\le \ell$; in particular
$1=x_1^0\dots x_{\ell}^0$ is a term. A term with an associated coefficient from $R$ is called \emph{monomial}. We endow the set of terms with the \emph{total degree-lexicographical order}, in which we first compare the degrees of the monomials and then break the ties by means of lexicographic order for the order $x_1>x_2>\dots>x_{\ell}$ on the variables. If the variables are only $x,y$ and $x>y,$ this gives the total order:
\begin{equation*}
1<y<x<y^2<xy<x^2<y^3<xy^2<x^2y<x^3<y^4<\cdots.
\end{equation*}
For $P\in R[\underline{X}]$ we denote by $HT(P)$ the leading term and by $HM(P)$ the highest monomial in $P$ with respect to the above order. We denote by $HC(P)$ the coefficient of the leading monomial in $P$. We denote by $T(P)$ the set of terms appearing in $P$ and by $M(P)$ the set of monomials in $P$.

For a given ideal $I\subseteq R[\underline{X}]$ an associated Gr\"obner basis is introduced as in Chapters 5, 10 in \cite{B}).

If $G\subseteq R[\underline{X}]$ is a finite set, we denote by $I(G)$ the ideal generated by $G$ in $R[\underline{X}]$.

\begin{definition} Let $f,g,p\in R[\underline{X}]$. We say that $f$ $D$-reduces to $g$ modulo $p$ and write $f \underset{p}{\to} g$ if there exists $m\in M(f)$ with $HM(p)\vert m$, say $m=m'\cdot HM(p),$ and $g=f-m'p$. For a finite set $G\subseteq R[\underline{X}]$, we denote by $\overset{*}{\underset{G}{\to}}$ the reflexive-transitive closure of $\underset{p}{\to}, p\in G$. We say that $g$ is a normal form for $f$ with respect to $G$ if $f \overset{*}{\underset{G}{\to}} g$ and no further $D$-reduction is possible. We say that $f$ is $D$-reducible modulo $G$ if $f\overset{*}{\underset{G}{\to}} 0$.
\end{definition}

If $f\overset{*}{\underset{G}{\to}} 0$, then $f\in I(G)$. The converse is also true if $G$ is a Gr\"obner basis.

\begin{definition} A $D$-Gr\"obner basis is a finite set $G$ of $R[\underline{X}]$ with the property that all $D$-normal forms modulo $G$ of elements of $I(G)$ equal zero. If $I\subseteq R[\underline{X}]$ is an ideal, then a $D$-Gr\"obner basis of $I$ is a $D$-Gr\"obner basis that generates the ideal $I$.
\end{definition}

\begin{propo} Let $G$ be a finite set of $R[\underline{X}]$. Then the following statements are equivalent:
\begin{enumerate}
\item $G$ is a Gr\"obner basis.
\item Every $f\not = 0, f\in I(G),$ is $D$-reducible modulo $G$.
\end{enumerate}
\end{propo}

We observe, nevertheless, that if $R$ is only a (PID), the normal form associated to a polynomial $f$ by a finite set $G\subseteq R[\underline{X}]$ is not unique. That is, the reminder of the division of $f$ by $G$ is not unique.

We introduce now the notions of $S$-polynomial and $G$-polynomial that allows to check if a given finite set $G\subseteq R[\underline{X}]$ is a Gr\"obner basis for the ideal it generates. As usual, $\emph{lcm}$ is the notation for the least common multiple and $\emph{gcd}$ is the notation for the greatest common divisor.

\begin{definition} Let $0\not = g_i\in R[\underline{X}], i=1,2,$ with $HC(g_i)=a_i$ and $HT(g_i)=t_i$. Let $a=b_ia_i=\text{lcm}(a_1,a_2)$ with $b_i\in R$, and $t=s_it_i=\text{lcm}(t_1,t_2)$ with $s_i\in T$. The the $S$-polynomial of $g_1,g_2$ is defined as:
\begin{equation*}
S(g_1,g_2)=b_1s_1g_1-b_2s_2g_2.
\end{equation*}

If $c_1,c_2\in R$ such that $\text{gcd}(a_1,a_2)=c_1a_1+c_2a_2$. Then the $G$-polynomial of $g_1,g_2$ is defined as:
\begin{equation*}
G(g_1,g_2)=c_1s_1g_1+c_2s_2g_2.
\end{equation*}
\end{definition}

\begin{thm} Let $G$ be a finite set of $R[\underline{X}]$. Assume that for all $g_1,g_2\in G$, $S(g_1,g_2)\overset{*}{\underset{G}{\to}} 0$ and $G(g_1,g_2)$ is top-$D$-reducible modulo $G$. Then $G$ is a Gr\"obner basis.
\end{thm}

Assume now that $R$ is an Euclidean domain with unique reminders (see~\cite[p. 463]{B}). This is the case for the ring of integers $\mathbb{Z}$ if we specify reminders upon division by $0\not =m$ to be in the interval $[0,m)$.

\begin{definition} Let $f,g,p\in R[\underline{X}]$. We say that $f$ $E$-reduces to $g$ modulo $p$ and write $f \underset{p}{\to} g$ if there exists $m=at\in M(f)$ with $HM(p)\vert t$, say $t=s\cdot HT(p),$ and $g=f-qsp$ where $0\not =q\in R$ is the quotient of $a$ upon division with unique reminder by $HC(p)$.
\end{definition}

\begin{propo} $E$-reduction extends $D$-reduction, i.e., every $D$-reduction step in an $E$-reduction step.
\end{propo}

\begin{thm} Let $R$ be an Euclidean domain with unique reminders, and assume $G\subseteq R[\underline{X}]$ is a $D$-Gr\"obner basis. Then the following hold:
\begin{enumerate}
\item $f \overset{*}{\underset{G}{\to}} 0$ for all $f\in I(G),$ where $\overset{*}{\underset{G}{\to}}$ denotes the $E$-reduction modulo $G$.
\item $E$-reduction modulo $G$ has unique normal forms.
\end{enumerate}
\end{thm}

The following result connect signed tilings and Gr\"obner bases. See~\cite{B-N} and~\cite{n-bone} for a proof.

\begin{thm}\label{thm-groebner-tiling} A polyomino $P$ admits a signed tiling by translates of prototiles
$P_1, P_2, \dots , P_k$ if and only if for some (test) monomial $x^{\alpha}y^{\beta}$ the polynomial
$x^{\alpha}y^{\beta}f_P$ is in the ideal generated in $\mathbb{Z}[X,Y]$ by $f_{P_1},\dots, f_{P_k}$.
\end{thm}

\section{Gr\"obner basis for $\mathcal{T}_n, n$ even}\label{s:3}

We show Gr\"obner bases for the ideals generated by $\mathcal{T}_4, \mathcal{T}_6$, as these are different from the general case.

\begin{propo} The  polynomials $C_1(2)=x^2+x+y+1,\  C_2(2)=y^2+x+y+1$ form a Gr\"obner basis for the ideal generated by $\mathcal{T}_4$.
\end{propo}

\begin{pf} The polynomials corresponding to tiles in $\mathcal{T}_4$ are $C_1(2), C_2(2), xy^2+xy+y^2+x, x^2y+xy+x^2+y.$
The last two can be generated by $C_1(2), C_2(2)$:
\begin{equation*}
\begin{gathered}
xy^2+xy+y^2+x=-C_1(2)+(x+1)C_2(2)\\
x^2y+xy+x^2+y=-C_2(2)+(y+1)C_1(2).
\end{gathered}
\end{equation*}

It remains to show that the $S$-polynomial associated to $C_1(2),C_2(2)$ can be reduced. One has:
\begin{equation*}
S(C_1(2),C_2(2))=y^2(x^2+x+y+1)-x^2(y^2+y+1+x)=(x+y+1)C_1(2)+(x+y+1)C_2(2).
\end{equation*}

\end{pf}

\begin{propo}  A Gr\"obner basis for the ideal of polynomials generated by $\mathcal{T}_6$ is given by:
\begin{equation*}
\begin{gathered}
C_1(3)=x^3+x^2+x+y^2+y+1,\ \ \ C_2(3)=y^3+y^2+y+x^2+x+1,\ \ \ C_3(3)=xy-1.
\end{gathered}
\end{equation*}
\end{propo}

\begin{pf} The polynomials associated to $\mathcal{T}_6$ are:
\begin{equation}\label{eq:poly-t6}
\begin{aligned}
H_1(k)&=y^4+y^3+y^2+y+1+x,\ \ H_2(k)=y^{4}+xy^4+xy^3+xy^2+xy+x,\\
H_3(k)&=y+x^4+x^3+x^2+x+1,\ \ H_4(k)=x^4y+x^3y+x^2y+xy+y+x^{4}.
\end{aligned}
\end{equation}

\begin{figure}[h!]
\centering
\begin{subfigure}{.20\textwidth}
\begin{tikzpicture}[scale=.38]

\draw [line width = 1,<->] (0,5.5)--(0,0)--(8,0);
\draw [line width = 2] (0,0)--(0,2)--(1,2)--(1,1)--(5,1)--(5,0)--(0,0);
\draw [line width = 1] (1,0)--(1,1);
\draw [line width = 1] (0,1)--(1,1);
\draw [line width = 1] (2,0)--(2,1);
\draw [line width = 1] (3,0)--(3,1);
\draw [line width = 1] (4,0)--(4,1);

\node at (.5,.5) {$-$};
\node at (1.5,.5) {$-$};
\node at (2.5,.5) {$-$};
\node at (3.5,.5) {$-$};
\node at (4.5,.5) {$-$};
\node at (.5,1.5) {$-$};
\end{tikzpicture}
\caption{Step 1 (-)}
\end{subfigure}
~~~~
\begin{subfigure}{.20\textwidth}
\begin{tikzpicture}[scale=.38]

\draw [line width = 1,<->] (0,5.5)--(0,0)--(8,0);
\draw [line width = 2] (4,0)--(4,1)--(0,1)--(0,2)--(5,2)--(5,0)--(4,0);
\draw [line width = 1] (0,0)--(0,2)--(1,2)--(1,1)--(5,1)--(5,0)--(0,0);
\draw [line width = 1] (1,0)--(1,1);
\draw [line width = 1] (0,1)--(1,1);
\draw [line width = 1] (2,0)--(2,1);
\draw [line width = 1] (3,0)--(3,1);
\draw [line width = 1] (4,0)--(4,1);

\draw [line width = 1] (4,1)--(4,2);
\draw [line width = 1] (2,1)--(2,2);
\draw [line width = 1] (3,1)--(3,2);
\draw [line width = 1] (4,1)--(5,1);

\node at (.5,.5) {$-$};
\node at (1.5,.5) {$-$};
\node at (2.5,.5) {$-$};
\node at (3.5,.5) {$-$};
\node at (4.5,.5) {$0$};
\node at (.5,1.5) {$0$};

\node at (4.5,1.5) {$+$};
\node at (1.5,1.5) {$+$};
\node at (2.5,1.5) {$+$};
\node at (3.5,1.5) {$+$};
\end{tikzpicture}
\caption{Step 2 (+)}
\end{subfigure}
~~~~
\begin{subfigure}{.20\textwidth}
\begin{tikzpicture}[scale=.38]

\draw [line width = 1,<->] (0,5.5)--(0,0)--(8,0);
\draw [line width = 2] (1,0)--(3,0)--(3,1)--(2,1)--(2,5)--(1,5)--(1,0);
\draw [line width = 1] (4,0)--(4,1)--(0,1)--(0,2)--(5,2)--(5,0)--(4,0);
\draw [line width = 1] (0,0)--(0,2)--(1,2)--(1,1)--(5,1)--(5,0)--(0,0);
\draw [line width = 1] (1,0)--(1,1);
\draw [line width = 1] (0,1)--(1,1);
\draw [line width = 1] (2,0)--(2,1);
\draw [line width = 1] (3,0)--(3,1);
\draw [line width = 1] (4,0)--(4,1);

\draw [line width = 1] (4,1)--(4,2);
\draw [line width = 1] (2,1)--(2,2);
\draw [line width = 1] (3,1)--(3,2);
\draw [line width = 1] (4,1)--(5,1);

\draw [line width = 1] (1,3)--(2,3);
\draw [line width = 1] (1,4)--(2,4);

\node at (1.5,2.5) {$+$};
\node at (1.5,3.5) {$+$};
\node at (1.5,4.5) {$+$};

\node at (.5,.5) {$-$};
\node at (1.5,.5) {$0$};
\node at (2.5,.5) {$0$};
\node at (3.5,.5) {$-$};
\node at (4.5,.5) {$0$};
\node at (.5,1.5) {$0$};

\node at (4.5,1.5) {$+$};
\node at (1.5,1.5) {\tiny{$+2$}};
\node at (2.5,1.5) {$+$};
\node at (3.5,1.5) {$+$};
\end{tikzpicture}
\caption{Step 3 (+)}
\end{subfigure}
~~~~
\begin{subfigure}{.20\textwidth}
\begin{tikzpicture}[scale=.38]

\draw [line width = 1,<->] (0,5.5)--(0,0)--(8,0);
\draw [line width = 2] (2,0)--(3,0)--(3,5)--(1,5)--(1,4)--(2,4)--(2,0);
\draw [line width = 1] (1,0)--(3,0)--(3,1)--(2,1)--(2,5)--(1,5)--(1,0);
\draw [line width = 1] (4,0)--(4,1)--(0,1)--(0,2)--(5,2)--(5,0)--(4,0);
\draw [line width = 1] (0,0)--(0,2)--(1,2)--(1,1)--(5,1)--(5,0)--(0,0);
\draw [line width = 1] (1,0)--(1,1);
\draw [line width = 1] (0,1)--(1,1);
\draw [line width = 1] (2,0)--(2,1);
\draw [line width = 1] (3,0)--(3,1);
\draw [line width = 1] (4,0)--(4,1);

\draw [line width = 1] (4,1)--(4,2);
\draw [line width = 1] (2,1)--(2,2);
\draw [line width = 1] (3,1)--(3,2);
\draw [line width = 1] (4,1)--(5,1);

\draw [line width = 1] (1,3)--(3,3);
\draw [line width = 1] (1,4)--(3,4);

\node at (1.5,2.5) {$+$};
\node at (1.5,3.5) {$+$};
\node at (1.5,4.5) {$0$};
\node at (2.5,4.5) {$-$};
\node at (2.5,3.5) {$-$};
\node at (2.5,2.5) {$-$};

\node at (.5,.5) {$-$};
\node at (1.5,.5) {$0$};
\node at (2.5,.5) {$-$};
\node at (3.5,.5) {$-$};
\node at (4.5,.5) {$0$};
\node at (.5,1.5) {$0$};

\node at (4.5,1.5) {$+$};
\node at (1.5,1.5) {\tiny{$+2$}};
\node at (2.5,1.5) {$0$};
\node at (3.5,1.5) {$+$};
\end{tikzpicture}
\caption{Step 4 (-)}
\end{subfigure}
~~~~
\begin{subfigure}{.20\textwidth}
\begin{tikzpicture}[scale=.38]

\draw [line width = 1,<->] (0,5.5)--(0,0)--(8,0);
\draw [line width = 2] (2,0)--(4,0)--(4,1)--(3,1)--(3,5)--(2,5)--(2,0);
\draw [line width = 1] (2,0)--(3,0)--(3,5)--(1,5)--(1,4)--(2,4)--(2,0);
\draw [line width = 1] (1,0)--(3,0)--(3,1)--(2,1)--(2,5)--(1,5)--(1,0);
\draw [line width = 1] (4,0)--(4,1)--(0,1)--(0,2)--(5,2)--(5,0)--(4,0);
\draw [line width = 1] (0,0)--(0,2)--(1,2)--(1,1)--(5,1)--(5,0)--(0,0);
\draw [line width = 1] (1,0)--(1,1);
\draw [line width = 1] (0,1)--(1,1);
\draw [line width = 1] (2,0)--(2,1);
\draw [line width = 1] (3,0)--(3,1);
\draw [line width = 1] (4,0)--(4,1);

\draw [line width = 1] (4,1)--(4,2);
\draw [line width = 1] (2,1)--(2,2);
\draw [line width = 1] (3,1)--(3,2);
\draw [line width = 1] (4,1)--(5,1);

\draw [line width = 1] (1,3)--(3,3);
\draw [line width = 1] (1,4)--(3,4);

\node at (1.5,2.5) {$+$};
\node at (1.5,3.5) {$+$};
\node at (1.5,4.5) {$0$};
\node at (2.5,4.5) {$0$};
\node at (2.5,3.5) {$0$};
\node at (2.5,2.5) {$0$};

\node at (.5,.5) {$-$};
\node at (1.5,.5) {$0$};
\node at (2.5,.5) {$0$};
\node at (3.5,.5) {$0$};
\node at (4.5,.5) {$0$};
\node at (.5,1.5) {$0$};

\node at (4.5,1.5) {$+$};
\node at (1.5,1.5) {\tiny{$+2$}};
\node at (2.5,1.5) {$+$};
\node at (3.5,1.5) {$+$};
\end{tikzpicture}
\caption{Step 5 (+)}
\end{subfigure}
~~~~
\begin{subfigure}{.20\textwidth}
\begin{tikzpicture}[scale=.38]

\draw [line width = 1,<->] (0,5.5)--(0,0)--(8,0);

\draw [line width = 2] (1,1)--(6,1)--(6,2)--(2,2)--(2,3)--(1,3)--(1,1);
\draw [line width = 1] (2,0)--(4,0)--(4,1)--(3,1)--(3,5)--(2,5)--(2,0);
\draw [line width = 1] (2,0)--(3,0)--(3,5)--(1,5)--(1,4)--(2,4)--(2,0);
\draw [line width = 1] (1,0)--(3,0)--(3,1)--(2,1)--(2,5)--(1,5)--(1,0);
\draw [line width = 1] (4,0)--(4,1)--(0,1)--(0,2)--(5,2)--(5,0)--(4,0);
\draw [line width = 1] (0,0)--(0,2)--(1,2)--(1,1)--(5,1)--(5,0)--(0,0);
\draw [line width = 1] (1,0)--(1,1);
\draw [line width = 1] (0,1)--(1,1);
\draw [line width = 1] (2,0)--(2,1);
\draw [line width = 1] (3,0)--(3,1);
\draw [line width = 1] (4,0)--(4,1);

\draw [line width = 1] (4,1)--(4,2);
\draw [line width = 1] (2,1)--(2,2);
\draw [line width = 1] (3,1)--(3,2);
\draw [line width = 1] (4,1)--(5,1);

\draw [line width = 1] (1,3)--(3,3);
\draw [line width = 1] (1,4)--(3,4);

\node at (1.5,2.5) {$0$};
\node at (1.5,3.5) {$+$};
\node at (1.5,4.5) {$0$};
\node at (2.5,4.5) {$0$};
\node at (2.5,3.5) {$0$};
\node at (2.5,2.5) {$0$};

\node at (.5,.5) {$-$};
\node at (1.5,.5) {$0$};
\node at (2.5,.5) {$0$};
\node at (3.5,.5) {$0$};
\node at (4.5,.5) {$0$};
\node at (.5,1.5) {$0$};

\node at (4.5,1.5) {$0$};
\node at (1.5,1.5) {$+$};
\node at (2.5,1.5) {$0$};
\node at (3.5,1.5) {$0$};

\node at (5.5,1.5) {$-$};
\end{tikzpicture}
\caption{Step 6 (-)}
\end{subfigure}
~~~~
\begin{subfigure}{.20\textwidth}
\begin{tikzpicture}[scale=.38]

\draw [line width = 1,<->] (0,5.5)--(0,0)--(8,0);

\draw [line width = 2] (1,2)--(1,3)--(6,3)--(6,1)--(5,1)--(5,2)--(1,2);
\draw [line width = 1] (1,1)--(6,1)--(6,2)--(2,2)--(2,3)--(1,3)--(1,1);
\draw [line width = 1] (2,0)--(4,0)--(4,1)--(3,1)--(3,5)--(2,5)--(2,0);
\draw [line width = 1] (2,0)--(3,0)--(3,5)--(1,5)--(1,4)--(2,4)--(2,0);
\draw [line width = 1] (1,0)--(3,0)--(3,1)--(2,1)--(2,5)--(1,5)--(1,0);
\draw [line width = 1] (4,0)--(4,1)--(0,1)--(0,2)--(5,2)--(5,0)--(4,0);
\draw [line width = 1] (0,0)--(0,2)--(1,2)--(1,1)--(5,1)--(5,0)--(0,0);
\draw [line width = 1] (1,0)--(1,1);
\draw [line width = 1] (0,1)--(1,1);
\draw [line width = 1] (2,0)--(2,1);
\draw [line width = 1] (3,0)--(3,1);
\draw [line width = 1] (4,0)--(4,1);

\draw [line width = 1] (4,1)--(4,2);
\draw [line width = 1] (2,1)--(2,2);
\draw [line width = 1] (3,1)--(3,2);
\draw [line width = 1] (4,1)--(5,1);

\draw [line width = 1] (1,3)--(3,3);
\draw [line width = 1] (1,4)--(3,4);

\draw [line width = 1] (4,2)--(4,3);
\draw [line width = 1] (5,2)--(5,3);

\node at (1.5,2.5) {$+$};
\node at (1.5,3.5) {$+$};
\node at (1.5,4.5) {$0$};
\node at (2.5,4.5) {$0$};
\node at (2.5,3.5) {$0$};
\node at (2.5,2.5) {$+$};

\node at (3.5,2.5) {$+$};
\node at (4.5,2.5) {$+$};
\node at (5.5,2.5) {$+$};

\node at (.5,.5) {$-$};
\node at (1.5,.5) {$0$};
\node at (2.5,.5) {$0$};
\node at (3.5,.5) {$0$};
\node at (4.5,.5) {$0$};
\node at (.5,1.5) {$0$};

\node at (4.5,1.5) {$0$};
\node at (1.5,1.5) {$+$};
\node at (2.5,1.5) {$0$};
\node at (3.5,1.5) {$0$};

\node at (5.5,1.5) {$0$};
\end{tikzpicture}
\caption{Step 7 (+)}
\end{subfigure}
~~~~
\begin{subfigure}{.20\textwidth}
\begin{tikzpicture}[scale=.38]

\draw [line width = 1,<->] (0,5.5)--(0,0)--(8,0);

\draw [fill=lightgray] (0,0)--(1,0)--(1,1)--(0,1)--(0,0);
\draw [fill=lightgray] (1,1)--(2,1)--(2,2)--(1,2)--(1,1);

\draw [line width = 2] (1,2)--(6,2)--(6,3)--(2,3)--(2,4)--(1,4)--(1,2);
\draw [line width = 1] (1,2)--(1,3)--(6,3)--(6,1)--(5,1)--(5,2)--(1,2);
\draw [line width = 1] (1,1)--(6,1)--(6,2)--(2,2)--(2,3)--(1,3)--(1,1);
\draw [line width = 1] (2,0)--(4,0)--(4,1)--(3,1)--(3,5)--(2,5)--(2,0);
\draw [line width = 1] (2,0)--(3,0)--(3,5)--(1,5)--(1,4)--(2,4)--(2,0);
\draw [line width = 1] (1,0)--(3,0)--(3,1)--(2,1)--(2,5)--(1,5)--(1,0);
\draw [line width = 1] (4,0)--(4,1)--(0,1)--(0,2)--(5,2)--(5,0)--(4,0);
\draw [line width = 1] (0,0)--(0,2)--(1,2)--(1,1)--(5,1)--(5,0)--(0,0);
\draw [line width = 1] (1,0)--(1,1);
\draw [line width = 1] (0,1)--(1,1);
\draw [line width = 1] (2,0)--(2,1);
\draw [line width = 1] (3,0)--(3,1);
\draw [line width = 1] (4,0)--(4,1);

\draw [line width = 1] (4,1)--(4,2);
\draw [line width = 1] (2,1)--(2,2);
\draw [line width = 1] (3,1)--(3,2);
\draw [line width = 1] (4,1)--(5,1);

\draw [line width = 1] (1,3)--(3,3);
\draw [line width = 1] (1,4)--(3,4);

\draw [line width = 1] (4,2)--(4,3);
\draw [line width = 1] (5,2)--(5,3);

\node at (1.5,2.5) {$0$};
\node at (1.5,3.5) {$0$};
\node at (1.5,4.5) {$0$};
\node at (2.5,4.5) {$0$};
\node at (2.5,3.5) {$0$};
\node at (2.5,2.5) {$0$};

\node at (3.5,2.5) {$0$};
\node at (4.5,2.5) {$0$};
\node at (5.5,2.5) {$0$};

\node at (.5,.5) {$-$};
\node at (1.5,.5) {$0$};
\node at (2.5,.5) {$0$};
\node at (3.5,.5) {$0$};
\node at (4.5,.5) {$0$};
\node at (.5,1.5) {$0$};

\node at (4.5,1.5) {$0$};
\node at (1.5,1.5) {$+$};
\node at (2.5,1.5) {$0$};
\node at (3.5,1.5) {$0$};

\node at (5.5,1.5) {$0$};
\end{tikzpicture}
\caption{Step 8 (-)}
\end{subfigure}
\caption{The polynomial $C_3(6)$ is generated by $\{H_1(6),H_2(6),H_3(6),H_4(6)\}$.}
\label{fig:special-form-356}
\end{figure}
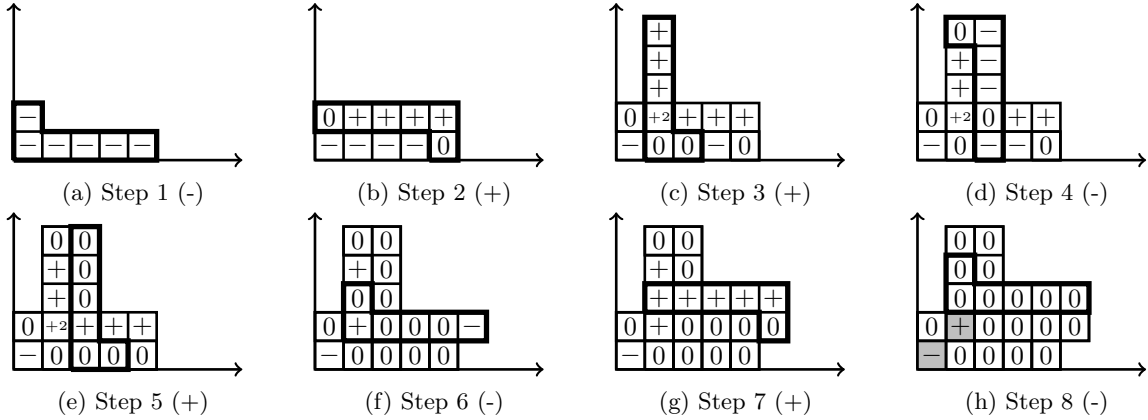

Similar to what is done in~\cite{Nit-odd}, the presence of $C_3(3)$ in the Gr\"obner basis allows to reduce the algebraic proofs to combinatorial considerations. We leave most of the details of this proof to the reader. The proof that \eqref{eq:poly-t6}  are in the ideal generated by $C_1(3),C_2(3),C_3(3)$  is similar to that of~\cite[Proposition 5]{Nit-odd}. The proof that $C_1(3),C_2(3),C_3(3)$ are in the ideal generated by  \eqref{eq:poly-t6}  is similar to that of~\cite[Proposition 6]{Nit-odd}. A geometric proof that $C_3(3)$ belongs to the ideal generated by $H_1(3), H_2(3), H_3(3),H_4(3)$ is shown in Figure~\ref{fig:special-form-356}.

\end{pf}

For the rest of this section $n=2k,$ where $k\ge 4$. The polynomials associated to the tiles in $\mathcal{T}_n$ are:

\begin{equation}\label{eq:generators-k-even}
\begin{aligned}
H_1(k)&=\frac{y^{2k-1}-1}{y-1}+x,\ H_2(k)=y^{2k-2}+\frac{x(y^{2k-1}-1)}{y-1},\\
H_3(k)&=y+\frac{x^{2k-1}-1}{x-1},\ H_4(k)=\frac{y(x^{2k-1}-1)}{x-1}+x^{2k-2}.
\end{aligned}
\end{equation}

We show that a Gr\"obner basis for the ideal generated in $\mathbb{Z}[X,Y]$ by \eqref{eq:generators-k-even} is given by:

\begin{equation}\label{eq:basis-k}
\begin{aligned}
C_1(k)&=\frac{y^{k+1}-1}{y-1}+x\cdot\frac{x^{k-1}-1}{x-1}+\left \lfloor{\frac{k-1}{2}}\right \rfloor xy-\left \lfloor{\frac{k-1}{2}}\right \rfloor,\\
C_2(k)&=\frac{x^{k+1}-1}{x-1}+y\cdot\frac{y^{k-1}-1}{y-1}+\left \lfloor{\frac{k-1}{2}}\right \rfloor xy-\left \lfloor{\frac{k-1}{2}}\right \rfloor,\\
C_3(k)&=x^2y + xy - x - 1,\ \ C_4(k)=xy^2 + xy - y - 1,\ \ C_5(k)=(k-2)xy-(k-2),
\end{aligned}
\end{equation}
where $\lfloor x \rfloor$ is the integer part of $x$.

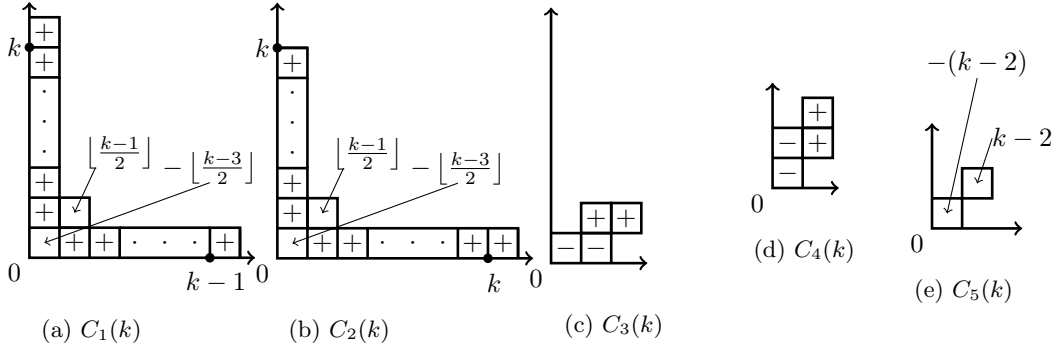
\begin{figure}[h!]
\centering
\begin{subfigure}{.15\textwidth}
\begin{tikzpicture}[scale=.4]
\draw [line width = 1,<->] (0,8.5)--(0,0)--(7.5,0);
\draw [line width = 1] (0,1)--(3,1)--(3,0);
\draw [line width = 1] (1,0)--(1,3)--(0,3);
\draw [line width = 1] (2,0)--(2,1);
\draw [line width = 1] (0,2)--(1,2);
\draw [line width = 1] (6,0)--(6,1)--(7,1)--(7,0);
\draw [line width = 1] (0,6)--(1,6)--(1,8)--(0,8);
\draw [line width = 1] (7,0)--(7,1);
\draw [line width = 1] (0,7)--(1,7);
\draw [line width = 1] (1,3)--(1,6);
\draw [line width = 1] (3,1)--(6,1);
\draw [line width = 1] (1,1)--(2,1)--(2,2)--(1,2)--(1,1);

\node at (.5,.5) {};
\node at (1.5,.5) {$+$};
\node at (2.5,.5) {$+$};
\node at (6.5,.5) {$+$};

\node at (.5,1.5) {$+$};
\node at (.5,2.5) {$+$};
\node at (.5,6.5) {$+$};
\node at (.5,7.5) {$+$};

\node at (3.5,.5) {$\cdot$};
\node at (4.5,.5) {$\cdot$};
\node at (5.5,.5) {$\cdot$};

\node at (.5,3.5) {$\cdot$};
\node at (.5,4.5) {$\cdot$};
\node at (.5,5.5) {$\cdot$};

\draw [fill] (6,0) circle (4pt);
\draw [fill] (0,7) circle (4pt);

\node at (-.5, -.5) {$0$};
\node at (6.2, -.8) {$k-1$};
\node at (-.5, 7) {$k$};
\node at (6,3) {$-\left \lfloor{\frac{k-3}{2}}\right \rfloor$};
\draw [->] (6,2.5)--(.5,.5);

\node at (3,3.5) {$\left \lfloor{\frac{k-1}{2}}\right \rfloor$};
\draw [<-] (1.5,1.5)--(2.2,3);
\end{tikzpicture}
\caption{$C_1(k)$}
\end{subfigure}
~~~~~
\begin{subfigure}{.15\textwidth}
\begin{tikzpicture}[scale=.4]
\draw [line width = 1,<->] (0,8.5)--(0,0)--(8.5,0);
\draw [line width = 1] (0,1)--(3,1)--(3,0);
\draw [line width = 1] (1,0)--(1,3)--(0,3);
\draw [line width = 1] (2,0)--(2,1);
\draw [line width = 1] (0,2)--(1,2);
\draw [line width = 1] (6,0)--(6,1)--(8,1)--(8,0);
\draw [line width = 1] (0,6)--(1,6)--(1,7)--(0,7);
\draw [line width = 1] (7,0)--(7,1);
\draw [line width = 1] (8,0)--(8,1);
\draw [line width = 1] (0,7)--(1,7);
\draw [line width = 1] (1,3)--(1,6);
\draw [line width = 1] (3,1)--(6,1);
\draw [line width = 1] (1,1)--(2,1)--(2,2)--(1,2)--(1,1);

\node at (.5,.5) {};
\node at (1.5,.5) {$+$};
\node at (2.5,.5) {$+$};
\node at (6.5,.5) {$+$};
\node at (7.5,.5) {$+$};

\node at (.5,1.5) {$+$};
\node at (.5,2.5) {$+$};
\node at (.5,6.5) {$+$};


\node at (3.5,.5) {$\cdot$};
\node at (4.5,.5) {$\cdot$};
\node at (5.5,.5) {$\cdot$};

\node at (.5,3.5) {$\cdot$};
\node at (.5,4.5) {$\cdot$};
\node at (.5,5.5) {$\cdot$};

\draw [fill] (7,0) circle (4pt);
\draw [fill] (0,7) circle (4pt);

\node at (-.5, -.5) {$0$};
\node at (7.2, -.8) {$k$};
\node at (-.5, 7) {$k$};
\node at (6,3) {$-\left \lfloor{\frac{k-3}{2}}\right \rfloor$};
\draw [->] (6,2.5)--(.5,.5);

\node at (3,3.5) {$\left \lfloor{\frac{k-1}{2}}\right \rfloor$};
\draw [<-] (1.5,1.5)--(2.2,3);
\end{tikzpicture}
\caption{$C_2(k)$}
\end{subfigure}
~~~~~~~~
\begin{subfigure}{.15\textwidth}
\begin{tikzpicture}[scale=.4]
\draw [line width = 1,<->] (0,8.5)--(0,0)--(3.2,0);
\draw [line width = 1] (0,1)--(1,1)--(1,0);
\draw [line width = 1] (1,1)--(1,2)--(2,2)--(2,1)--(1,1);
\draw [line width = 1] (1,1)--(2,1)--(2,0);
\draw [line width = 1] (2,1)--(2,2)--(3,2)--(3,1)--(2,1);

\node at (.5,.5) {$-$};
\node at (1.5,1.5) {$+$};
\node at (1.5,.5) {$-$};
\node at (2.5,1.5) {$+$};
\node at (-.5, -.5) {$0$};
\end{tikzpicture}
\caption{$C_3(k)$}
\end{subfigure}
~~
\begin{subfigure}{.1\textwidth}
\begin{tikzpicture}[scale=.4]
\draw [line width = 1,<->] (0,3.5)--(0,0)--(2.2,0);
\draw [line width = 1] (0,1)--(1,1)--(1,0);
\draw [line width = 1] (1,1)--(1,2)--(2,2)--(2,1)--(1,1);
\draw [line width = 1] (0,2)--(1,2)--(1,1);
\draw [line width = 1] (1,2)--(1,3)--(2,3)--(2,2)--(1,2);

\node at (.5,.5) {$-$};
\node at (1.5,1.5) {$+$};
\node at (.5,1.5) {$-$};
\node at (1.5,2.5) {$+$};
\node at (-.5, -.5) {$0$};
\end{tikzpicture}
\caption{$C_4(k)$}
\end{subfigure}
~~
\begin{subfigure}{.1\textwidth}
\begin{tikzpicture}[scale=.4]
\draw [line width = 1,<->] (0,3.5)--(0,0)--(3,0);
\draw [line width = 1] (0,1)--(1,1)--(1,0);
\draw [line width = 1] (1,1)--(1,2)--(2,2)--(2,1)--(1,1);

\node at (.5,.5) {};
\node at (1.5,1.5) {};
\node at (-.5, -.5) {$0$};

\node at (1.5,5.5) {$-(k-2)$};
\draw [->] (1.5,5)--(.5,.5);

\node at (3,3) {$k-2$};
\draw [->] (2,3)--(1.5,1.5);
\end{tikzpicture}
\caption{$C_5(k)$}
\end{subfigure}
\caption{The Gr\"obner basis $\{C_1(k), C_2(k), C_3(k),C_4(k),C_5(k)\}.$}
\label{fig:signed-tiles-even}
\end{figure}

It is convenient to visualize the elements of the basis as tiles with cells labeled by integers, see Figure~\ref{fig:signed-tiles-even}.

\begin{propo}\label{p:prop5-even} The polynomials $H_i(k), 1\le i\le 4,$ belong to the ideal generated by $C_i(k), 1\le i\le 5$.
\end{propo}

\begin{pf} Due to the symmetry, it is enough to show that $H_1(k), H_2(k)$ belong to the ideal. The polynomials $C_3(k), C_4(k)$ allow to translate an horizontal domino with both two cells labeled by the same sign, respectively a vertical domino, along a vector parallel to the first bisector $y=x$. They also allow to translate horizontally or vertically a block of two cells adjacent at a vertex and labeled by different signs into a similar block. If the length of the translation is even, the signs stay the same. If the length of the translation is odd, all signs are changed. See Figure~\ref{fig:diagonaltiles}.

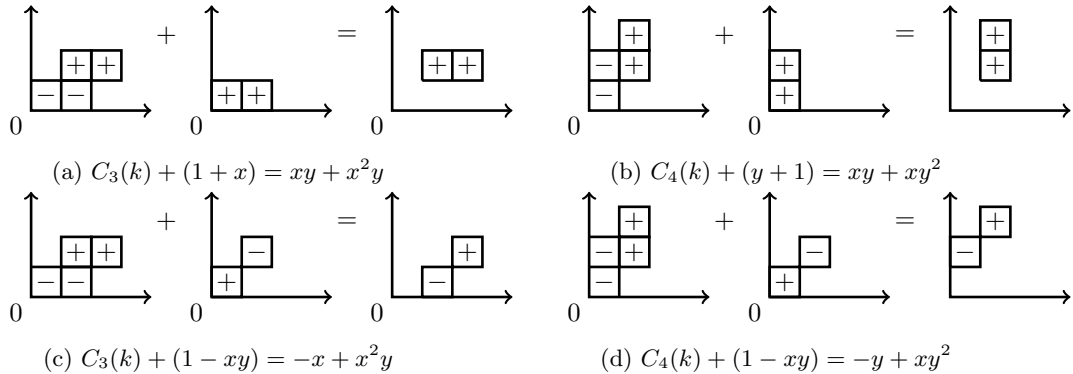
\begin{figure}[h!]
\centering
\begin{subfigure}{.35\textwidth}
\begin{tikzpicture}[scale=.4]
\draw [line width = 1,<->] (0,3.5)--(0,0)--(4,0);
\draw [line width = 1] (0,1)--(1,1)--(1,0);
\draw [line width = 1] (1,1)--(1,2)--(2,2)--(2,1)--(1,1);
\draw [line width = 1] (1,1)--(2,1)--(2,0);
\draw [line width = 1] (2,1)--(2,2)--(3,2)--(3,1)--(2,1);

\draw [line width = 1,<->] (6,3.5)--(6,0)--(10,0);
\draw [line width = 1] (6,0)--(8,0)--(8,1)--(6,1)--(6,0);
\draw [line width = 1] (7,0)--(7,1);

\node at (6.5,.5) {$+$};
\node at (7.5,.5) {$+$};
\node at (5.5, -.5) {$0$};

\node at (4.5,2.5) {$+$};
\node at (10.5,2.5) {$=$};

\draw [line width = 1,<->] (12,3.5)--(12,0)--(16,0);
\draw [line width = 1] (13,1)--(15,1)--(15,2)--(13,2)--(13,1);
\draw [line width = 1] (14,1)--(14,2);

\node at (13.5,1.5) {$+$};
\node at (14.5,1.5) {$+$};
\node at (11.5, -.5) {$0$};

\node at (.5,.5) {$-$};
\node at (1.5,1.5) {$+$};
\node at (1.5,.5) {$-$};
\node at (2.5,1.5) {$+$};
\node at (-.5, -.5) {$0$};
\end{tikzpicture}
\caption{$C_3(k)+(1+x)=xy+x^2y$}
\end{subfigure}
~~~~~~~~~~~~
\begin{subfigure}{.35\textwidth}
\begin{tikzpicture}[scale=.4]
\draw [line width = 1,<->] (0,3.5)--(0,0)--(4,0);
\draw [line width = 1] (0,1)--(1,1)--(1,0);
\draw [line width = 1] (1,1)--(1,2)--(2,2)--(2,1)--(1,1);
\draw [line width = 1] (0,2)--(1,2)--(1,1);
\draw [line width = 1] (1,2)--(1,3)--(2,3)--(2,2)--(1,2);

\draw [line width = 1,<->] (6,3.5)--(6,0)--(10,0);
\draw [line width = 1] (6,0)--(7,0)--(7,2)--(6,2)--(6,0);
\draw [line width = 1] (6,1)--(7,1);

\draw [line width = 1,<->] (12,3.5)--(12,0)--(16,0);
\draw [line width = 1] (13,1)--(14,1)--(14,3)--(13,3)--(13,1);
\draw [line width = 1] (13,2)--(14,2);

\node at (13.5,1.5) {$+$};
\node at (13.5,2.5) {$+$};

\node at (6.5,.5) {$+$};
\node at (6.5,1.5) {$+$};
\node at (5.5, -.5) {$0$};

\node at (4.5,2.5) {$+$};
\node at (10.5,2.5) {$=$};

\node at (.5,.5) {$-$};
\node at (1.5,1.5) {$+$};
\node at (.5,1.5) {$-$};
\node at (1.5,2.5) {$+$};
\node at (-.5, -.5) {$0$};
\end{tikzpicture}
\caption{$C_4(k)+(y+1)=xy+xy^2$}
\end{subfigure}
~~~~~~~~~~~
\begin{subfigure}{.35\textwidth}
\begin{tikzpicture}[scale=.4]
\draw [line width = 1,<->] (0,3.5)--(0,0)--(4,0);
\draw [line width = 1] (0,1)--(1,1)--(1,0);
\draw [line width = 1] (1,1)--(1,2)--(2,2)--(2,1)--(1,1);
\draw [line width = 1] (1,1)--(2,1)--(2,0);
\draw [line width = 1] (2,1)--(2,2)--(3,2)--(3,1)--(2,1);

\draw [line width = 1,<->] (6,3.5)--(6,0)--(10,0);
\draw [line width = 1] (6,0)--(7,0)--(7,2)--(8,2)--(8,1)--(6,1)--(6,0);

\node at (6.5,.5) {$+$};
\node at (7.5,1.5) {$-$};
\node at (5.5, -.5) {$0$};

\node at (4.5,2.5) {$+$};
\node at (10.5,2.5) {$=$};

\draw [line width = 1,<->] (12,3.5)--(12,0)--(16,0);
\draw [line width = 1] (13,0)--(14,0)--(14,2)--(15,2)--(15,1)--(13,1)--(13,0);

\node at (13.5,.5) {$-$};
\node at (14.5,1.5) {$+$};
\node at (11.5, -.5) {$0$};

\node at (.5,.5) {$-$};
\node at (1.5,1.5) {$+$};
\node at (1.5,.5) {$-$};
\node at (2.5,1.5) {$+$};
\node at (-.5, -.5) {$0$};
\end{tikzpicture}
\caption{$C_3(k)+(1-xy)=-x+x^2y$}
\end{subfigure}
~~~~~~~~~~~~
\begin{subfigure}{.35\textwidth}
\begin{tikzpicture}[scale=.4]
\draw [line width = 1,<->] (0,3.5)--(0,0)--(4,0);
\draw [line width = 1] (0,1)--(1,1)--(1,0);
\draw [line width = 1] (1,1)--(1,2)--(2,2)--(2,1)--(1,1);
\draw [line width = 1] (0,2)--(1,2)--(1,1);
\draw [line width = 1] (1,2)--(1,3)--(2,3)--(2,2)--(1,2);

\draw [line width = 1,<->] (6,3.5)--(6,0)--(10,0);
\draw [line width = 1] (6,0)--(7,0)--(7,2)--(8,2)--(8,1)--(6,1)--(6,0);

\draw [line width = 1,<->] (12,3.5)--(12,0)--(16,0);
\draw [line width = 1] (12,1)--(13,1)--(13,3)--(14,3)--(14,2)--(12,2)--(12,1);

\node at (12.5,1.5) {$-$};
\node at (13.5,2.5) {$+$};

\node at (6.5,.5) {$+$};
\node at (7.5,1.5) {$-$};
\node at (5.5, -.5) {$0$};

\node at (4.5,2.5) {$+$};
\node at (10.5,2.5) {$=$};

\node at (.5,.5) {$-$};
\node at (1.5,1.5) {$+$};
\node at (.5,1.5) {$-$};
\node at (1.5,2.5) {$+$};
\node at (-.5, -.5) {$0$};
\end{tikzpicture}
\caption{$C_4(k)+(1-xy)=-y+xy^2$}
\end{subfigure}
\caption{Tiles arithmetic.}
\label{fig:diagonaltiles}
\end{figure}

We show how to build $H_1(k)$. There are two cases to be considered, $k$ odd and $k$ even.

The steps of a geometric constructions for $k$ odd are shown in Figure~\ref{fig:H1-new(k)kodd}. To reach Step 1, we add several times multiples of $C_4(k)$, as in Figure~\ref{fig:diagonaltiles}, b). To reach Step 2, we add several times multiples of $C_3(k)$, as in Figure~\ref{fig:diagonaltiles}, a). To reach Step 3, first we subtract $C_5(k)$, then add several times multiples of $C_3(k), C_4(k)$ as in Figure~\ref{fig:diagonaltiles}, c), d). To obtain now $H_1(k)$ in the initial position, we multiply the tile in Step 3 by $x^{k-2}$, which will translate the tile $k-2$ cells up, and then add multiples on $C_3(k), C_4(k)$, as in Figure~\ref{fig:diagonaltiles}, c), d).

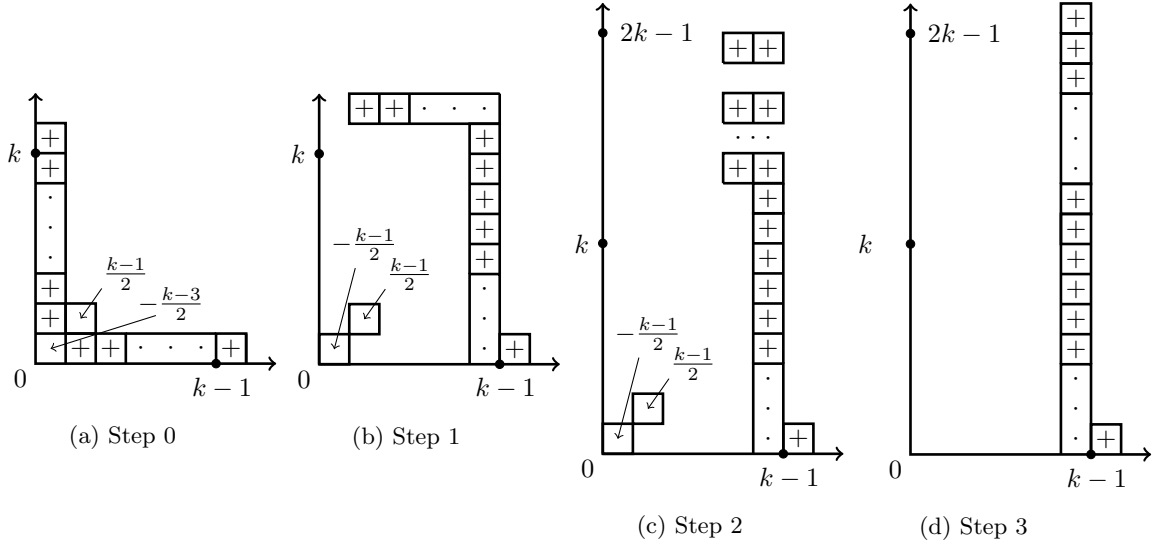
\begin{figure}[h!]
\centering
\begin{subfigure}{.2\textwidth}
\begin{tikzpicture}[scale=.4]
\draw [line width = 1,<->] (0,9)--(0,0)--(8,0);
\draw [line width = 1] (0,1)--(3,1)--(3,0);
\draw [line width = 1] (1,0)--(1,3)--(0,3);
\draw [line width = 1] (2,0)--(2,1);
\draw [line width = 1] (0,2)--(1,2);
\draw [line width = 1] (6,0)--(6,1)--(7,1)--(7,0);
\draw [line width = 1] (0,6)--(1,6)--(1,8)--(0,8);
\draw [line width = 1] (7,0)--(7,1);
\draw [line width = 1] (0,7)--(1,7);
\draw [line width = 1] (1,3)--(1,6);
\draw [line width = 1] (3,1)--(6,1);
\draw [line width = 1] (1,1)--(2,1)--(2,2)--(1,2)--(1,1);

\node at (.5,.5) {};
\node at (1.5,.5) {$+$};
\node at (2.5,.5) {$+$};
\node at (6.5,.5) {$+$};

\node at (.5,1.5) {$+$};
\node at (.5,2.5) {$+$};
\node at (.5,6.5) {$+$};
\node at (.5,7.5) {$+$};

\node at (3.5,.5) {$\cdot$};
\node at (4.5,.5) {$\cdot$};
\node at (5.5,.5) {$\cdot$};

\node at (.5,3.5) {$\cdot$};
\node at (.5,4.5) {$\cdot$};
\node at (.5,5.5) {$\cdot$};

\draw [fill] (6,0) circle (4pt);
\draw [fill] (0,7) circle (4pt);

\node at (-.5, -.5) {$0$};
\node at (6.2, -.8) {$k-1$};
\node at (-.7, 7) {$k$};
\node at (4.5,2) {$-\frac{k-3}{2}$};
\draw [->] (4,1.6)--(.5,.5);

\node at (3,3) {$\frac{k-1}{2}$};
\draw [<-] (1.5,1.5)--(2.2,2.5);
\end{tikzpicture}
\caption{Step 0}
\end{subfigure}
~~
\begin{subfigure}{.2\textwidth}
\begin{tikzpicture}[scale=.4]
\draw [line width = 1,<->] (0,9)--(0,0)--(8,0);

\draw [line width = 1] (6,0)--(6,1)--(7,1)--(7,0);

\draw [line width = 1] (1,1)--(2,1)--(2,2)--(1,2)--(1,1);
\draw [line width = 1] (0,0)--(1,0)--(1,1)--(0,1)--(0,0);

\draw [line width = 1] (5,0)--(5,3)--(6,3)--(6,0);
\draw [line width = 1] (5,3)--(5,8)--(6,8)--(6,3);
\draw [line width = 1] (5,4)--(6,4);
\draw [line width = 1] (5,5)--(6,5);
\draw [line width = 1] (5,6)--(6,6);
\draw [line width = 1] (5,7)--(6,7);

\draw [line width = 1] (6,9)--(1,9)--(1,8)--(6,8)--(6,9);
\draw [line width = 1] (2,8)--(2,9);
\draw [line width = 1] (3,8)--(3,9);

\node at (.5,.5) {};
\node at (1.5,8.5) {$+$};
\node at (2.5,8.5) {$+$};
\node at (5.5,7.5) {$+$};
\node at (5.5,6.5) {$+$};
\node at (5.5,5.5) {$+$};
\node at (6.5,.5) {$+$};
\node at (5.5,4.5) {$+$};
\node at (5.5,3.5) {$+$};

\node at (5.5,.5) {$\cdot$};
\node at (5.5,1.5) {$\cdot$};
\node at (5.5,2.5) {$\cdot$};

\node at (3.5,8.5) {$\cdot$};
\node at (4.5,8.5) {$\cdot$};
\node at (5.5,8.5) {$\cdot$};

\draw [fill] (6,0) circle (4pt);
\draw [fill] (0,7) circle (4pt);

\node at (-.5, -.5) {$0$};
\node at (6.2, -.8) {$k-1$};
\node at (-.7, 7) {$k$};
\node at (1.5,4) {$-\frac{k-1}{2}$};
\draw [->] (1.2,3.2)--(.5,.5);

\node at (3,3) {$\frac{k-1}{2}$};
\draw [<-] (1.5,1.5)--(2.2,2.5);
\end{tikzpicture}
\caption{Step 1}
\end{subfigure}
~~
\begin{subfigure}{.2\textwidth}
\begin{tikzpicture}[scale=.4]
\draw [line width = 1,<->] (0,15)--(0,0)--(8,0);

\draw [line width = 1] (6,0)--(6,1)--(7,1)--(7,0);

\draw [line width = 1] (1,1)--(2,1)--(2,2)--(1,2)--(1,1);
\draw [line width = 1] (0,0)--(1,0)--(1,1)--(0,1)--(0,0);

\draw [line width = 1] (5,0)--(5,3)--(6,3)--(6,0);
\draw [line width = 1] (5,3)--(5,8)--(6,8)--(6,3);
\draw [line width = 1] (5,4)--(6,4);
\draw [line width = 1] (5,5)--(6,5);
\draw [line width = 1] (5,6)--(6,6);
\draw [line width = 1] (5,7)--(6,7);

\draw [line width = 1] (6,9)--(5,9)--(5,8)--(6,8)--(6,9);

\draw [line width = 1] (4,9)--(6,9)--(6,10)--(4,10)--(4,9);
\draw [line width = 1] (5,9)--(5,10);

\node at (4.5,9.5) {$+$};
\node at (5.5,9.5) {$+$};

\draw [line width = 1] (4,11)--(6,11)--(6,12)--(4,12)--(4,11);
\draw [line width = 1] (5,11)--(5,12);

\node at (4.5,11.5) {$+$};
\node at (5.5,11.5) {$+$};

\draw [line width = 1] (4,13)--(6,13)--(6,14)--(4,14)--(4,13);
\draw [line width = 1] (5,13)--(5,14);

\node at (4.5,13.5) {$+$};
\node at (5.5,13.5) {$+$};

\node at (4.5,10.5) {$\cdot$};
\node at (5,10.5) {$\cdot$};
\node at (5.5,10.5) {$\cdot$};

\node at (5.5,7.5) {$+$};
\node at (5.5,6.5) {$+$};
\node at (5.5,5.5) {$+$};
\node at (6.5,.5) {$+$};
\node at (5.5,4.5) {$+$};
\node at (5.5,3.5) {$+$};

\node at (5.5,.5) {$\cdot$};
\node at (5.5,1.5) {$\cdot$};
\node at (5.5,2.5) {$\cdot$};

\node at (5.5,8.5) {$+$};

\draw [fill] (6,0) circle (4pt);
\draw [fill] (0,7) circle (4pt);
\draw [fill] (0,14) circle (4pt);

\node at (-.5, -.5) {$0$};
\node at (6.2, -.8) {$k-1$};
\node at (-.7, 7) {$k$};
\node at (1.5,4) {$-\frac{k-1}{2}$};
\draw [->] (1.2,3.2)--(.5,.5);
\node at (1.8, 14) {$2k-1$};

\node at (3,3) {$\frac{k-1}{2}$};
\draw [<-] (1.5,1.5)--(2.2,2.5);
\end{tikzpicture}
\caption{Step 2}
\end{subfigure}
~~
\begin{subfigure}{.2\textwidth}
\begin{tikzpicture}[scale=.4]
\draw [line width = 1,<->] (0,15)--(0,0)--(8,0);

\draw [line width = 1] (6,0)--(6,1)--(7,1)--(7,0);

\draw [line width = 1] (5,0)--(5,3)--(6,3)--(6,0);
\draw [line width = 1] (5,3)--(5,8)--(6,8)--(6,3);
\draw [line width = 1] (5,4)--(6,4);
\draw [line width = 1] (5,5)--(6,5);
\draw [line width = 1] (5,6)--(6,6);

\draw [line width = 1] (5,8)--(5,7)--(6,7);

\draw [line width = 1] (5,12)--(6,12)--(6,13)--(5,13)--(5,12);

\node at (5.5,12.5) {$+$};

\draw [line width = 1] (5,14)--(6,14)--(6,15)--(5,15)--(5,14);
\draw [line width = 1] (5,13)--(6,13)--(6,14)--(5,14)--(5,13);
\draw [line width = 1] (5,8)--(6,8)--(6,9)--(5,9)--(5,8);

\draw [line width = 1] (5,9)--(6,9)--(6,12)--(5,12)--(5,9);

\node at (5.5,14.5) {$+$};
\node at (5.5,13.5) {$+$};

\node at (5.5,9.5) {$\cdot$};
\node at (5.5,11.5) {$\cdot$};
\node at (5.5,10.5) {$\cdot$};

\node at (5.5,7.5) {$+$};
\node at (5.5,6.5) {$+$};
\node at (5.5,5.5) {$+$};
\node at (6.5,.5) {$+$};
\node at (5.5,4.5) {$+$};
\node at (5.5,3.5) {$+$};

\node at (5.5,.5) {$\cdot$};
\node at (5.5,1.5) {$\cdot$};
\node at (5.5,2.5) {$\cdot$};

\node at (5.5,8.5) {$+$};

\draw [fill] (6,0) circle (4pt);
\draw [fill] (0,7) circle (4pt);
\draw [fill] (0,14) circle (4pt);

\node at (-.5, -.5) {$0$};
\node at (6.2, -.8) {$k-1$};
\node at (-1.5, 7) {$k$};

\node at (1.8, 14) {$2k-1$};
\end{tikzpicture}
\caption{Step 3}
\end{subfigure}
\caption{Building $H_1(k),k$ odd, out of $\{C_1(k), C_2(k), C_3(k),C_4(k),C_5(k)\}.$}
\label{fig:H1-new(k)kodd}
\end{figure}

The steps of a geometric constructions for $k$ even are shown in Figure~\ref{fig:H1-new(k)keven}. To reach Step 1, we add several times multiples of $C_4(k)$, as in Figure~\ref{fig:diagonaltiles}, b). To reach Step 2, we add several times multiples of $C_3(k)$, as in Figure~\ref{fig:diagonaltiles}, a). To reach Step 3, first we subtract $C_5(k)$, then add several times multiples of $C_3(k), C_4(k)$ as in Figure~\ref{fig:diagonaltiles}, c), d). To obtain now $H_1(k)$ in the initial position, we multiply the tile in Step 3 by $x^{k-2}$, which will translate the tile $k-2$ cells up, and then add multiples on $C_3(k), C_4(k)$, as in Figure~\ref{fig:diagonaltiles}, c), d).

\begin{figure}[h!]
\centering
\begin{subfigure}{.2\textwidth}
\begin{tikzpicture}[scale=.4]
\draw [line width = 1,<->] (0,9)--(0,0)--(8,0);
\draw [line width = 1] (0,1)--(3,1)--(3,0);
\draw [line width = 1] (1,0)--(1,3)--(0,3);
\draw [line width = 1] (2,0)--(2,1);
\draw [line width = 1] (0,2)--(1,2);
\draw [line width = 1] (6,0)--(6,1)--(7,1)--(7,0);
\draw [line width = 1] (0,6)--(1,6)--(1,8)--(0,8);
\draw [line width = 1] (7,0)--(7,1);
\draw [line width = 1] (0,7)--(1,7);
\draw [line width = 1] (1,3)--(1,6);
\draw [line width = 1] (3,1)--(6,1);
\draw [line width = 1] (1,1)--(2,1)--(2,2)--(1,2)--(1,1);

\node at (.5,.5) {};
\node at (1.5,.5) {$+$};
\node at (2.5,.5) {$+$};
\node at (6.5,.5) {$+$};

\node at (.5,1.5) {$+$};
\node at (.5,2.5) {$+$};
\node at (.5,6.5) {$+$};
\node at (.5,7.5) {$+$};

\node at (3.5,.5) {$\cdot$};
\node at (4.5,.5) {$\cdot$};
\node at (5.5,.5) {$\cdot$};

\node at (.5,3.5) {$\cdot$};
\node at (.5,4.5) {$\cdot$};
\node at (.5,5.5) {$\cdot$};

\draw [fill] (6,0) circle (4pt);
\draw [fill] (0,7) circle (4pt);

\node at (-.5, -.5) {$0$};
\node at (6.2, -.8) {$k-1$};
\node at (-.7, 7) {$k$};
\node at (4.5,2) {$-\frac{k-3}{2}$};
\draw [->] (4,1.6)--(.5,.5);

\node at (3,3) {$\frac{k-1}{2}$};
\draw [<-] (1.5,1.5)--(2.2,2.5);
\end{tikzpicture}
\caption{Step 0}
\end{subfigure}
~~
\begin{subfigure}{.2\textwidth}
\begin{tikzpicture}[scale=.4]
\draw [line width = 1,<->] (0,9)--(0,0)--(8,0);

\draw [line width = 1] (6,0)--(6,1)--(7,1)--(7,0);

\draw [line width = 1] (1,1)--(2,1)--(2,2)--(1,2)--(1,1);
\draw [line width = 1] (0,0)--(1,0)--(1,1)--(0,1)--(0,0);

\draw [line width = 1] (5,0)--(5,3)--(6,3)--(6,0);
\draw [line width = 1] (5,3)--(5,8)--(6,8)--(6,3);
\draw [line width = 1] (5,4)--(6,4);
\draw [line width = 1] (5,5)--(6,5);
\draw [line width = 1] (5,6)--(6,6);

\draw [line width = 1] (5,8)--(0,8)--(0,7)--(6,7);
\draw [line width = 1] (1,7)--(1,8);
\draw [line width = 1] (2,7)--(2,8);

\node at (.5,.5) {};
\node at (.5,7.5) {$+$};
\node at (1.5,7.5) {$+$};
\node at (5.5,7.5) {$+$};
\node at (5.5,6.5) {$+$};
\node at (5.5,5.5) {$+$};
\node at (6.5,.5) {$+$};
\node at (5.5,4.5) {$+$};
\node at (5.5,3.5) {$+$};

\node at (5.5,1.5) {$\cdot$};
\node at (5.5,2.5) {$\cdot$};
\node at (5.5,.5) {$\cdot$};

\node at (3.5,7.5) {$\cdot$};
\node at (4.5,7.5) {$\cdot$};
\node at (2.5,7.5) {$\cdot$};

\draw [fill] (6,0) circle (4pt);
\draw [fill] (0,7) circle (4pt);

\node at (-.5, -.5) {$0$};
\node at (6.2, -.8) {$k-1$};
\node at (-.7, 7) {$k$};
\node at (1.5,4) {$-\frac{k-1}{2}$};
\draw [->] (1.2,3.2)--(.5,.5);

\node at (3,3) {$\frac{k-1}{2}$};
\draw [<-] (1.5,1.5)--(2.2,2.5);
\end{tikzpicture}
\caption{Step 1}
\end{subfigure}
~~
\begin{subfigure}{.2\textwidth}
\begin{tikzpicture}[scale=.4]
\draw [line width = 1,<->] (0,15)--(0,0)--(8,0);

\draw [line width = 1] (6,0)--(6,1)--(7,1)--(7,0);

\draw [line width = 1] (1,1)--(2,1)--(2,2)--(1,2)--(1,1);
\draw [line width = 1] (0,0)--(1,0)--(1,1)--(0,1)--(0,0);

\draw [line width = 1] (5,0)--(5,3)--(6,3)--(6,0);
\draw [line width = 1] (5,3)--(5,8)--(6,8)--(6,3);
\draw [line width = 1] (5,4)--(6,4);
\draw [line width = 1] (5,5)--(6,5);
\draw [line width = 1] (5,6)--(6,6);

\draw [line width = 1] (5,8)--(4,8)--(4,7)--(6,7);
\draw [line width = 1] (4,9)--(6,9)--(6,10)--(4,10)--(4,9);
\draw [line width = 1] (5,9)--(5,10);

\node at (4.5,9.5) {$+$};
\node at (5.5,9.5) {$+$};

\draw [line width = 1] (4,11)--(6,11)--(6,12)--(4,12)--(4,11);
\draw [line width = 1] (5,11)--(5,12);

\node at (4.5,11.5) {$+$};
\node at (5.5,11.5) {$+$};

\draw [line width = 1] (4,13)--(6,13)--(6,14)--(4,14)--(4,13);
\draw [line width = 1] (5,13)--(5,14);

\node at (4.5,13.5) {$+$};
\node at (5.5,13.5) {$+$};

\node at (4.5,10.5) {$\cdot$};
\node at (5,10.5) {$\cdot$};
\node at (5.5,10.5) {$\cdot$};

\node at (5.5,7.5) {$+$};
\node at (5.5,6.5) {$+$};
\node at (5.5,5.5) {$+$};
\node at (6.5,.5) {$+$};
\node at (5.5,4.5) {$+$};
\node at (5.5,3.5) {$+$};

\node at (5.5,1.5) {$\cdot$};
\node at (5.5,2.5) {$\cdot$};
\node at (5.5,.5) {$\cdot$};

\node at (4.5,7.5) {$+$};

\draw [fill] (6,0) circle (4pt);
\draw [fill] (0,7) circle (4pt);
\draw [fill] (0,14) circle (4pt);

\node at (-.5, -.5) {$0$};
\node at (6.2, -.8) {$k-1$};
\node at (-.7, 7) {$k$};
\node at (1.5,4) {$-\frac{k-1}{2}$};
\draw [->] (1.2,3.2)--(.5,.5);
\node at (1.8, 14) {$2k-1$};

\node at (3,3) {$\frac{k-1}{2}$};
\draw [<-] (1.5,1.5)--(2.2,2.5);
\end{tikzpicture}
\caption{Step 2}
\end{subfigure}
~~
\begin{subfigure}{.2\textwidth}
\begin{tikzpicture}[scale=.4]
\draw [line width = 1,<->] (0,15)--(0,0)--(8,0);

\draw [line width = 1] (6,0)--(6,1)--(7,1)--(7,0);

\draw [line width = 1] (5,0)--(5,3)--(6,3)--(6,0);
\draw [line width = 1] (5,3)--(5,8)--(6,8)--(6,3);
\draw [line width = 1] (5,4)--(6,4);
\draw [line width = 1] (5,5)--(6,5);
\draw [line width = 1] (5,6)--(6,6);

\draw [line width = 1] (5,8)--(5,7)--(6,7);

\draw [line width = 1] (5,12)--(6,12)--(6,13)--(5,13)--(5,12);

\node at (5.5,12.5) {$+$};

\draw [line width = 1] (5,14)--(6,14)--(6,15)--(5,15)--(5,14);
\draw [line width = 1] (5,13)--(6,13)--(6,14)--(5,14)--(5,13);
\draw [line width = 1] (5,8)--(6,8)--(6,9)--(5,9)--(5,8);

\draw [line width = 1] (5,9)--(6,9)--(6,12)--(5,12)--(5,9);

\node at (5.5,14.5) {$+$};
\node at (5.5,13.5) {$+$};

\node at (5.5,9.5) {$\cdot$};
\node at (5.5,11.5) {$\cdot$};
\node at (5.5,10.5) {$\cdot$};

\node at (5.5,7.5) {$+$};
\node at (5.5,6.5) {$+$};
\node at (5.5,5.5) {$+$};
\node at (6.5,.5) {$+$};
\node at (5.5,4.5) {$+$};
\node at (5.5,3.5) {$+$};

\node at (5.5,.5) {$\cdot$};
\node at (5.5,1.5) {$\cdot$};
\node at (5.5,2.5) {$\cdot$};

\node at (5.5,8.5) {$+$};

\draw [fill] (6,0) circle (4pt);
\draw [fill] (0,7) circle (4pt);
\draw [fill] (0,14) circle (4pt);

\node at (-.5, -.5) {$0$};
\node at (6.2, -.8) {$k-1$};
\node at (-.7, 7) {$k$};

\node at (1.8, 14) {$2k-1$};
\end{tikzpicture}
\caption{Step 3}
\end{subfigure}
\caption{Building $H_1(k),k$ even, out of $\{C_1(k), C_2(k), C_3(k),C_4(k),C_5(k)\}.$}
\label{fig:H1-new(k)keven}
\end{figure}
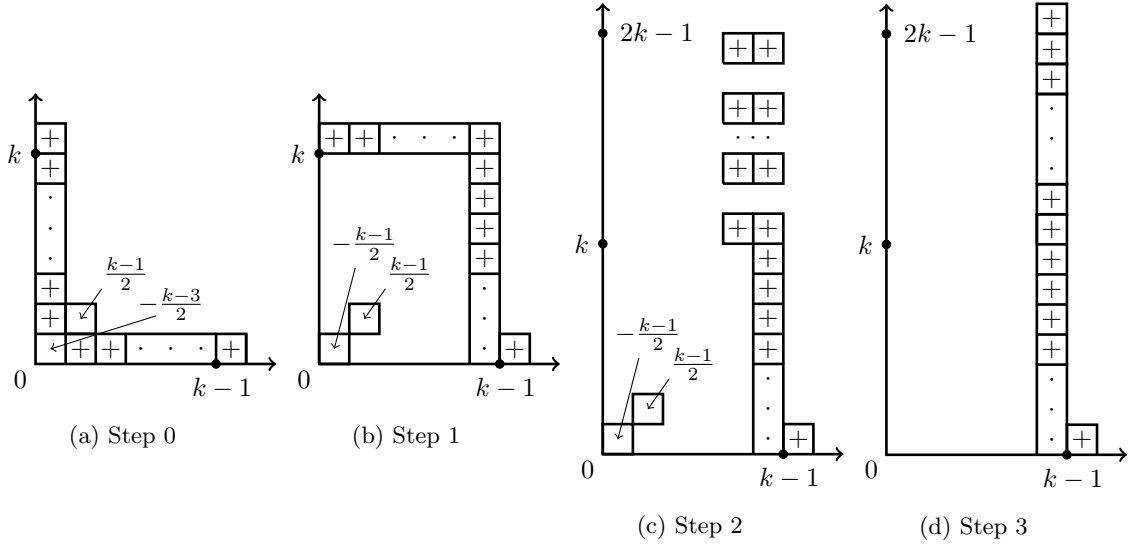

We show how to build $H_2(k)$. There are two cases to be considered, $k$ odd and $k$ even.

The steps of a geometric constructions for $k$ odd are shown in Figure~\ref{fig:H1(k)kodd}. To reach Step 1, we add several times multiples of $C_4(k)$, as in Figure~\ref{fig:diagonaltiles}, b). To reach Step 2, we add several times multiples of $C_3(k)$, as in Figure~\ref{fig:diagonaltiles}, a). To reach Step 3, first we subtract $C_5(k)$, then add several times multiples of $C_3(k), C_4(k)$ as in Figure~\ref{fig:diagonaltiles}, c), d). To obtain now $H_2(k)$ in the initial position, we multiply the tile in Step 3 by $x^{k-2}$, which will translate the tile $k-2$ cells up, and then add multiples on $C_3(k), C_4(k)$, as in Figure~\ref{fig:diagonaltiles}, c), d).

\begin{figure}[h!]
\centering
\begin{subfigure}{.2\textwidth}
\begin{tikzpicture}[scale=.4]
\draw [line width = 1,<->] (0,9)--(0,0)--(8,0);
\draw [line width = 1] (0,1)--(3,1)--(3,0);
\draw [line width = 1] (1,0)--(1,3)--(0,3);
\draw [line width = 1] (2,0)--(2,1);
\draw [line width = 1] (0,2)--(1,2);
\draw [line width = 1] (6,0)--(6,1)--(7,1)--(7,0);
\draw [line width = 1] (0,6)--(1,6)--(1,8)--(0,8);
\draw [line width = 1] (7,0)--(7,1);
\draw [line width = 1] (0,7)--(1,7);
\draw [line width = 1] (1,3)--(1,6);
\draw [line width = 1] (3,1)--(6,1);
\draw [line width = 1] (1,1)--(2,1)--(2,2)--(1,2)--(1,1);

\node at (.5,.5) {};
\node at (1.5,.5) {$+$};
\node at (2.5,.5) {$+$};
\node at (6.5,.5) {$+$};

\node at (.5,1.5) {$+$};
\node at (.5,2.5) {$+$};
\node at (.5,6.5) {$+$};
\node at (.5,7.5) {$+$};

\node at (3.5,.5) {$\cdot$};
\node at (4.5,.5) {$\cdot$};
\node at (5.5,.5) {$\cdot$};

\node at (.5,3.5) {$\cdot$};
\node at (.5,4.5) {$\cdot$};
\node at (.5,5.5) {$\cdot$};

\draw [fill] (6,0) circle (4pt);
\draw [fill] (0,7) circle (4pt);

\node at (-.5, -.5) {$0$};
\node at (6.2, -.8) {$k-1$};
\node at (-.7, 7) {$k$};
\node at (4.5,2) {$-\frac{k-4}{2}$};
\draw [->] (4,1.6)--(.5,.5);

\node at (3,3) {$\frac{k-2}{2}$};
\draw [<-] (1.5,1.5)--(2.2,2.5);
\end{tikzpicture}
\caption{Step 0}
\end{subfigure}
~~
\begin{subfigure}{.2\textwidth}
\begin{tikzpicture}[scale=.4]
\draw [line width = 1,<->] (0,9)--(0,0)--(8,0);

\draw [line width = 1] (6,0)--(6,1)--(7,1)--(7,0);

\draw [line width = 1] (1,1)--(2,1)--(2,2)--(1,2)--(1,1);
\draw [line width = 1] (0,0)--(1,0)--(1,1)--(0,1)--(0,0);

\draw [line width = 1] (6,1)--(6,4)--(7,4)--(7,1);
\draw [line width = 1] (6,4)--(6,9)--(7,9)--(7,4);
\draw [line width = 1] (6,5)--(7,5);
\draw [line width = 1] (6,6)--(7,6);
\draw [line width = 1] (6,7)--(7,7);

\draw [line width = 1] (6,9)--(1,9)--(1,8)--(7,8);
\draw [line width = 1] (2,8)--(2,9);
\draw [line width = 1] (3,8)--(3,9);

\node at (.5,.5) {};
\node at (1.5,8.5) {$+$};
\node at (2.5,8.5) {$+$};
\node at (6.5,8.5) {$+$};
\node at (6.5,7.5) {$+$};
\node at (6.5,6.5) {$+$};
\node at (6.5,.5) {$+$};
\node at (6.5,5.5) {$+$};
\node at (6.5,4.5) {$+$};

\node at (6.5,1.5) {$\cdot$};
\node at (6.5,2.5) {$\cdot$};
\node at (6.5,3.5) {$\cdot$};

\node at (3.5,8.5) {$\cdot$};
\node at (4.5,8.5) {$\cdot$};
\node at (5.5,8.5) {$\cdot$};

\draw [fill] (6,0) circle (4pt);
\draw [fill] (0,7) circle (4pt);

\node at (-.5, -.5) {$0$};
\node at (6.2, -.8) {$k-1$};
\node at (-.7, 7) {$k$};
\node at (1.5,4) {$-\frac{k-2}{2}$};
\draw [->] (1.2,3.2)--(.5,.5);

\node at (3,3) {$\frac{k-2}{2}$};
\draw [<-] (1.5,1.5)--(2.2,2.5);
\end{tikzpicture}
\caption{Step 1}
\end{subfigure}
~~
\begin{subfigure}{.2\textwidth}
\begin{tikzpicture}[scale=.4]
\draw [line width = 1,<->] (0,15)--(0,0)--(8,0);

\draw [line width = 1] (6,0)--(6,1)--(7,1)--(7,0);

\draw [line width = 1] (1,1)--(2,1)--(2,2)--(1,2)--(1,1);
\draw [line width = 1] (0,0)--(1,0)--(1,1)--(0,1)--(0,0);

\draw [line width = 1] (6,1)--(6,4)--(7,4)--(7,1);
\draw [line width = 1] (6,4)--(6,9)--(7,9)--(7,4);
\draw [line width = 1] (6,5)--(7,5);
\draw [line width = 1] (6,6)--(7,6);
\draw [line width = 1] (6,7)--(7,7);

\draw [line width = 1] (6,9)--(5,9)--(5,8)--(7,8);
\draw [line width = 1] (5,10)--(7,10)--(7,11)--(5,11)--(5,10);
\draw [line width = 1] (6,10)--(6,11);

\node at (5.5,10.5) {$+$};
\node at (6.5,10.5) {$+$};

\draw [line width = 1] (5,12)--(7,12)--(7,13)--(5,13)--(5,12);
\draw [line width = 1] (6,12)--(6,13);

\node at (5.5,12.5) {$+$};
\node at (6.5,12.5) {$+$};

\draw [line width = 1] (5,14)--(7,14)--(7,15)--(5,15)--(5,14);
\draw [line width = 1] (6,14)--(6,15);

\node at (5.5,14.5) {$+$};
\node at (6.5,14.5) {$+$};

\node at (5.5,11.5) {$\cdot$};
\node at (6,11.5) {$\cdot$};
\node at (6.5,11.5) {$\cdot$};

\node at (6.5,8.5) {$+$};
\node at (6.5,7.5) {$+$};
\node at (6.5,6.5) {$+$};
\node at (6.5,.5) {$+$};
\node at (6.5,5.5) {$+$};
\node at (6.5,4.5) {$+$};

\node at (6.5,1.5) {$\cdot$};
\node at (6.5,2.5) {$\cdot$};
\node at (6.5,3.5) {$\cdot$};

\node at (5.5,8.5) {$+$};

\draw [fill] (6,0) circle (4pt);
\draw [fill] (0,7) circle (4pt);
\draw [fill] (0,14) circle (4pt);

\node at (-.5, -.5) {$0$};
\node at (6.2, -.8) {$k-1$};
\node at (-.7, 7) {$k$};
\node at (1.5,4) {$-\frac{k-2}{2}$};
\draw [->] (1.2,3.2)--(.5,.5);
\node at (1.8, 14) {$2k-1$};

\node at (3,3) {$\frac{k-2}{2}$};
\draw [<-] (1.5,1.5)--(2.2,2.5);
\end{tikzpicture}
\caption{Step 2}
\end{subfigure}
~~
\begin{subfigure}{.2\textwidth}
\begin{tikzpicture}[scale=.4]
\draw [line width = 1,<->] (0,15)--(0,0)--(8,0);

\draw [line width = 1] (6,0)--(6,1)--(7,1)--(7,0);


\draw [line width = 1] (6,1)--(6,4)--(7,4)--(7,1);
\draw [line width = 1] (6,4)--(6,9)--(7,9)--(7,4);
\draw [line width = 1] (6,5)--(7,5);
\draw [line width = 1] (6,6)--(7,6);
\draw [line width = 1] (6,7)--(7,7);

\draw [line width = 1] (6,9)--(6,8)--(7,8);

\draw [line width = 1] (6,13)--(7,13)--(7,14)--(6,14)--(6,13);

\node at (6.5,13.5) {$+$};

\draw [line width = 1] (5,14)--(7,14)--(7,15)--(5,15)--(5,14);
\draw [line width = 1] (6,14)--(6,15);
\draw [line width = 1] (6,9)--(7,9)--(7,10)--(6,10)--(6,9);

\draw [line width = 1] (6,10)--(7,10)--(7,13)--(6,13)--(6,10);

\node at (5.5,14.5) {$+$};
\node at (6.5,14.5) {$+$};

\node at (6.5,12.5) {$\cdot$};
\node at (6.5,11.5) {$\cdot$};
\node at (6.5,10.5) {$\cdot$};

\node at (6.5,8.5) {$+$};
\node at (6.5,7.5) {$+$};
\node at (6.5,6.5) {$+$};
\node at (6.5,.5) {$+$};
\node at (6.5,5.5) {$+$};
\node at (6.5,4.5) {$+$};

\node at (6.5,1.5) {$\cdot$};
\node at (6.5,2.5) {$\cdot$};
\node at (6.5,3.5) {$\cdot$};

\node at (6.5,9.5) {$+$};

\draw [fill] (6,0) circle (4pt);
\draw [fill] (0,7) circle (4pt);
\draw [fill] (0,14) circle (4pt);

\node at (-.5, -.5) {$0$};
\node at (6.2, -.8) {$k-1$};
\node at (-.7, 7) {$k$};

\node at (1.8, 14) {$2k-1$};
\end{tikzpicture}
\caption{Step 3}
\end{subfigure}
\caption{Building $H_2(k),k$ odd, out of $\{C_1(k), C_2(k), C_3(k),C_4(k),C_5(k)\}.$}
\label{fig:H1(k)kodd}
\end{figure}
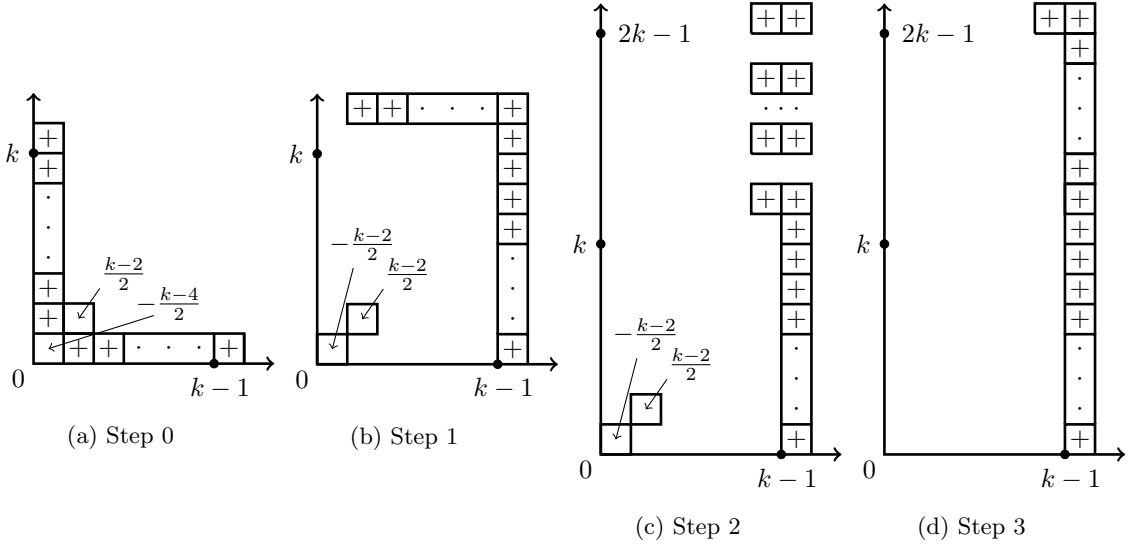

The steps of a geometric constructions for $k$ even are shown in Figure~\ref{fig:H1(k)keven}. To reach Step 1, we add several times multiples of $C_4(k)$, as in Figure~\ref{fig:diagonaltiles}, b). To reach Step 2, we add several times multiples of $C_3(k)$, as in Figure~\ref{fig:diagonaltiles}, a). To reach Step 3, first we subtract $C_5(k)$, then add several times multiples of $C_3(k), C_4(k)$ as in Figure~\ref{fig:diagonaltiles}, c), d). To obtain now $H_2(k)$ in the initial position, we multiply the tile in Step 3 by $x^{k-2}$, which will translate the tile $k-2$ cells up, and then add multiples on $C_3(k), C_4(k)$, as in Figure~\ref{fig:diagonaltiles}, c), d).

\begin{figure}[h!]
\centering
\begin{subfigure}{.2\textwidth}
\begin{tikzpicture}[scale=.4]
\draw [line width = 1,<->] (0,9)--(0,0)--(8,0);
\draw [line width = 1] (0,1)--(3,1)--(3,0);
\draw [line width = 1] (1,0)--(1,3)--(0,3);
\draw [line width = 1] (2,0)--(2,1);
\draw [line width = 1] (0,2)--(1,2);
\draw [line width = 1] (6,0)--(6,1)--(7,1)--(7,0);
\draw [line width = 1] (0,6)--(1,6)--(1,8)--(0,8);
\draw [line width = 1] (7,0)--(7,1);
\draw [line width = 1] (0,7)--(1,7);
\draw [line width = 1] (1,3)--(1,6);
\draw [line width = 1] (3,1)--(6,1);
\draw [line width = 1] (1,1)--(2,1)--(2,2)--(1,2)--(1,1);

\node at (.5,.5) {};
\node at (1.5,.5) {$+$};
\node at (2.5,.5) {$+$};
\node at (6.5,.5) {$+$};

\node at (.5,1.5) {$+$};
\node at (.5,2.5) {$+$};
\node at (.5,6.5) {$+$};
\node at (.5,7.5) {$+$};

\node at (3.5,.5) {$\cdot$};
\node at (4.5,.5) {$\cdot$};
\node at (5.5,.5) {$\cdot$};

\node at (.5,3.5) {$\cdot$};
\node at (.5,4.5) {$\cdot$};
\node at (.5,5.5) {$\cdot$};

\draw [fill] (6,0) circle (4pt);
\draw [fill] (0,7) circle (4pt);

\node at (-.5, -.5) {$0$};
\node at (6.2, -.8) {$k-1$};
\node at (-1.5, 7) {$k$};
\node at (4.5,2) {$-\frac{k-4}{2}$};
\draw [->] (4,1.6)--(.5,.5);

\node at (3,3) {$\frac{k-2}{2}$};
\draw [<-] (1.5,1.5)--(2.2,2.5);
\end{tikzpicture}
\caption{Step 0}
\end{subfigure}
~~~~~~~~~~
\begin{subfigure}{.25\textwidth}
\begin{tikzpicture}[scale=.4]
\draw [line width = 1,<->] (0,9)--(0,0)--(8,0);

\draw [line width = 1] (6,0)--(6,1)--(7,1)--(7,0);

\draw [line width = 1] (1,1)--(2,1)--(2,2)--(1,2)--(1,1);
\draw [line width = 1] (0,0)--(1,0)--(1,1)--(0,1)--(0,0);

\draw [line width = 1] (6,1)--(6,4)--(7,4)--(7,1);
\draw [line width = 1] (6,4)--(6,8)--(7,8)--(7,4);
\draw [line width = 1] (6,5)--(7,5);
\draw [line width = 1] (6,6)--(7,6);
\draw [line width = 1] (6,7)--(7,7);

\draw [line width = 1] (0,7)--(7,7)--(7,8)--(0,8)--(0,7);
\draw [line width = 1] (1,7)--(1,8);
\draw [line width = 1] (2,7)--(2,8);
\draw [line width = 1] (5,7)--(5,8);

\node at (.5,7.5) {$+$};
\node at (1.5,7.5) {$+$};
\node at (5.5,7.5) {$+$};
\node at (6.5,7.5) {$+$};
\node at (6.5,6.5) {$+$};
\node at (6.5,.5) {$+$};
\node at (6.5,5.5) {$+$};
\node at (6.5,4.5) {$+$};

\node at (6.5,1.5) {$\cdot$};
\node at (6.5,2.5) {$\cdot$};
\node at (6.5,3.5) {$\cdot$};

\node at (3.5,7.5) {$\cdot$};
\node at (4.5,7.5) {$\cdot$};
\node at (2.5,7.5) {$\cdot$};

\draw [fill] (6,0) circle (4pt);
\draw [fill] (0,7) circle (4pt);

\node at (-.5, -.5) {$0$};
\node at (6.2, -.8) {$k-1$};
\node at (-.7, 7) {$k$};
\node at (1.5,4) {$-\frac{k-2}{2}$};
\draw [->] (1.2,3.2)--(.5,.5);

\node at (3,3) {$\frac{k-2}{2}$};
\draw [<-] (1.5,1.5)--(2.2,2.5);
\end{tikzpicture}
\caption{Step 1}
\end{subfigure}
~~
\begin{subfigure}{.2\textwidth}
\begin{tikzpicture}[scale=.4]
\draw [line width = 1,<->] (0,15)--(0,0)--(8,0);

\draw [line width = 1] (6,0)--(6,1)--(7,1)--(7,0);

\draw [line width = 1] (1,1)--(2,1)--(2,2)--(1,2)--(1,1);
\draw [line width = 1] (0,0)--(1,0)--(1,1)--(0,1)--(0,0);

\draw [line width = 1] (6,1)--(6,4)--(7,4)--(7,1);
\draw [line width = 1] (6,4)--(6,9)--(7,9)--(7,4);
\draw [line width = 1] (6,5)--(7,5);
\draw [line width = 1] (6,6)--(7,6);
\draw [line width = 1] (6,7)--(7,7);

\draw [line width = 1] (6,9)--(5,9)--(5,8)--(7,8);
\draw [line width = 1] (5,10)--(7,10)--(7,11)--(5,11)--(5,10);
\draw [line width = 1] (6,10)--(6,11);

\node at (5.5,10.5) {$+$};
\node at (6.5,10.5) {$+$};

\draw [line width = 1] (5,12)--(7,12)--(7,13)--(5,13)--(5,12);
\draw [line width = 1] (6,12)--(6,13);

\node at (5.5,12.5) {$+$};
\node at (6.5,12.5) {$+$};

\draw [line width = 1] (5,14)--(7,14)--(7,15)--(5,15)--(5,14);
\draw [line width = 1] (6,14)--(6,15);

\node at (5.5,14.5) {$+$};
\node at (6.5,14.5) {$+$};

\node at (5.5,11.5) {$\cdot$};
\node at (6,11.5) {$\cdot$};
\node at (6.5,11.5) {$\cdot$};

\node at (6.5,8.5) {$+$};
\node at (6.5,7.5) {$+$};
\node at (6.5,6.5) {$+$};
\node at (6.5,.5) {$+$};
\node at (6.5,5.5) {$+$};
\node at (6.5,4.5) {$+$};

\node at (6.5,1.5) {$\cdot$};
\node at (6.5,2.5) {$\cdot$};
\node at (6.5,3.5) {$\cdot$};

\node at (5.5,8.5) {$+$};

\draw [fill] (6,0) circle (4pt);
\draw [fill] (0,8) circle (4pt);
\draw [fill] (0,14) circle (4pt);

\node at (-.5, -.5) {$0$};
\node at (6.2, -.8) {$k-1$};
\node at (-.7, 8) {$k$};
\node at (1.5,4) {$-\frac{k-2}{2}$};
\draw [->] (1.2,3.2)--(.5,.5);
\node at (1.8, 14) {$2k-1$};

\node at (3,3) {$\frac{k-2}{2}$};
\draw [<-] (1.5,1.5)--(2.2,2.5);
\end{tikzpicture}
\caption{Step 2}
\end{subfigure}
~~
\begin{subfigure}{.2\textwidth}
\begin{tikzpicture}[scale=.4]
\draw [line width = 1,<->] (0,15)--(0,0)--(8,0);

\draw [line width = 1] (6,0)--(6,1)--(7,1)--(7,0);


\draw [line width = 1] (6,1)--(6,4)--(7,4)--(7,1);
\draw [line width = 1] (6,4)--(6,9)--(7,9)--(7,4);
\draw [line width = 1] (6,5)--(7,5);
\draw [line width = 1] (6,6)--(7,6);
\draw [line width = 1] (6,7)--(7,7);

\draw [line width = 1] (6,9)--(6,8)--(7,8);

\draw [line width = 1] (6,13)--(7,13)--(7,14)--(6,14)--(6,13);

\node at (6.5,13.5) {$+$};

\draw [line width = 1] (5,14)--(7,14)--(7,15)--(5,15)--(5,14);
\draw [line width = 1] (6,14)--(6,15);
\draw [line width = 1] (6,9)--(7,9)--(7,10)--(6,10)--(6,9);

\draw [line width = 1] (6,10)--(7,10)--(7,13)--(6,13)--(6,10);

\node at (5.5,14.5) {$+$};
\node at (6.5,14.5) {$+$};

\node at (6.5,12.5) {$\cdot$};
\node at (6.5,11.5) {$\cdot$};
\node at (6.5,10.5) {$\cdot$};

\node at (6.5,8.5) {$+$};
\node at (6.5,7.5) {$+$};
\node at (6.5,6.5) {$+$};
\node at (6.5,.5) {$+$};
\node at (6.5,5.5) {$+$};
\node at (6.5,4.5) {$+$};

\node at (6.5,1.5) {$\cdot$};
\node at (6.5,2.5) {$\cdot$};
\node at (6.5,3.5) {$\cdot$};

\node at (6.5,9.5) {$+$};

\draw [fill] (6,0) circle (4pt);
\draw [fill] (0,7) circle (4pt);
\draw [fill] (0,14) circle (4pt);

\node at (-.5, -.5) {$0$};
\node at (6.2, -.8) {$k-1$};
\node at (-.7, 7) {$k$};

\node at (1.8, 14) {$2k-1$};
\end{tikzpicture}
\caption{Step 3}
\end{subfigure}
\caption{Building $H_2(k),k$ even, out of $\{C_1(k), C_2(k), C_3(k),C_4(k),C_5(k)\}.$}
\label{fig:H1(k)keven}
\end{figure}
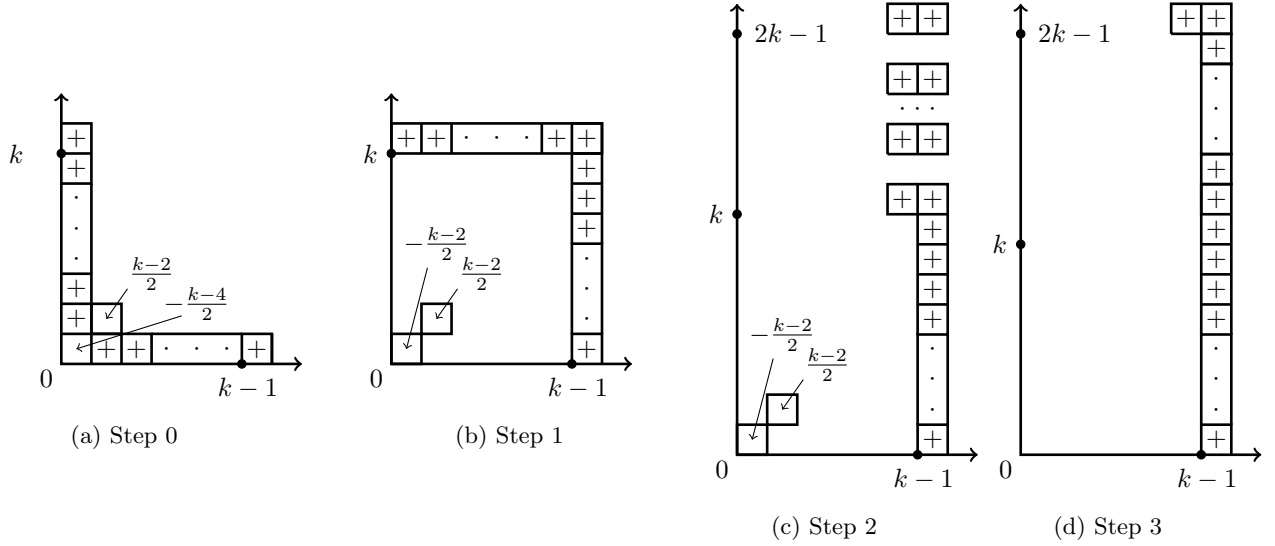
\end{pf}

\begin{propo}\label{p:prop6-even} The polynomials $C_i(k), 1\le i\le 5$ belong to the ideal generated by $H_i(k), 1\le i\le 4$.
\end{propo}

\begin{pf} Due to the symmetry, it is enough to show that $C_1(k), C_3(k)$ and $C_5(k)$ belong to the ideal. We show how to generate $C_3(k), C_5(k)$ (and consequently $C_4(k)$). To generate $C_1(k)$ we can reverse the process in Proposition~\ref{p:prop5-even}.
For $C_3(k)$, one has $C_3(k)=(xy+x-1)H_3(k)-xH_4(k).$ 
To generate $C_5(k)$ we first show how to obtain a configuration in which all nontrivial cells, 4 of them, are located on the main diagonal. See Figure~\ref{fig:main-diag345}. Then we use the tiles arithmetic shown in Figure~\ref{fig:diagonaltiles-bis22} to pull the cells in position $(k-1,k-1)$ and $(2k-2,2k-2)$ in positions $(1,1)$ and $(2,2)$. This constructs $C_5(k)$.

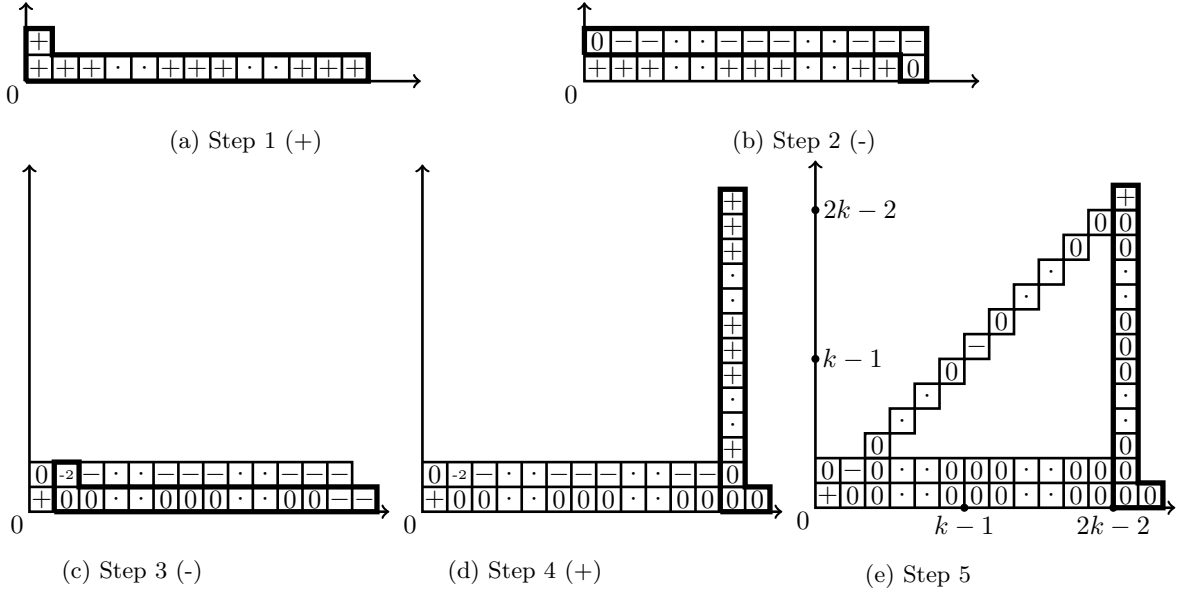
\begin{figure}[h!]
\centering
\begin{subfigure}{.4\textwidth}
\begin{tikzpicture}[scale=.35]
\draw [line width = 1,<->] (0,3)--(0,0)--(15,0);

\draw [line width = 2] (0,0)--(13,0)--(13,1)--(1,1)--(1,2)--(0,2)--(0,0);
\draw [line width = 1] (0,1)--(1,1);
\draw [line width = 1] (1,0)--(1,1);
\draw [line width = 1] (2,0)--(2,1);
\draw [line width = 1] (3,0)--(3,1);

\draw [line width = 1] (5,0)--(5,1);
\draw [line width = 1] (6,0)--(6,1);
\draw [line width = 1] (7,0)--(7,1);
\draw [line width = 1] (8,0)--(8,1);
\draw [line width = 1] (4,0)--(4,1);
\draw [line width = 1] (9,0)--(9,1);

\draw [line width = 1] (12,0)--(12,1);
\draw [line width = 1] (11,0)--(11,1);
\draw [line width = 1] (10,0)--(10,1);

\node at (.5,.5) {$+$};
\node at (1.5,.5) {$+$};
\node at (2.5,.5) {$+$};
\node at (10.5,.5) {$+$};
\node at (11.5,.5) {$+$};
\node at (12.5,.5) {$+$};

\node at (5.5,.5) {$+$};
\node at (6.5,.5) {$+$};
\node at (7.5,.5) {$+$};

\node at (.5,1.5) {$+$};

\node at (3.5,.5) {$\cdot$};
\node at (4.5,.5) {$\cdot$};

\node at (8.5,.5) {$\cdot$};
\node at (9.5,.5) {$\cdot$};

\node at (-.5, -.5) {$0$};
\end{tikzpicture}
\caption{Step 1 (+)}
\end{subfigure}
~~~~~
\begin{subfigure}{.4\textwidth}
\begin{tikzpicture}[scale=.35]
\draw [line width = 1,<->] (0,3)--(0,0)--(15,0);

\draw [line width = 2] (0,1)--(12,1)--(12,0)--(13,0)--(13,2)--(0,2)--(0,1);
\draw [line width = 1] (0,1)--(1,1);
\draw [line width = 1] (1,0)--(1,1);
\draw [line width = 1] (2,0)--(2,1);
\draw [line width = 1] (3,0)--(3,1);

\draw [line width = 1] (1,2)--(1,1);
\draw [line width = 1] (2,2)--(2,1);
\draw [line width = 1] (3,2)--(3,1);
\draw [line width = 1] (12,1)--(13,1);

\draw [line width = 1] (5,2)--(5,1);
\draw [line width = 1] (6,2)--(6,1);
\draw [line width = 1] (7,2)--(7,1);
\draw [line width = 1] (8,2)--(8,1);
\draw [line width = 1] (4,2)--(4,1);
\draw [line width = 1] (9,2)--(9,1);

\draw [line width = 1] (5,0)--(5,1);
\draw [line width = 1] (6,0)--(6,1);
\draw [line width = 1] (7,0)--(7,1);
\draw [line width = 1] (8,0)--(8,1);
\draw [line width = 1] (4,0)--(4,1);
\draw [line width = 1] (9,0)--(9,1);

\draw [line width = 1] (12,0)--(12,1);
\draw [line width = 1] (11,0)--(11,1);
\draw [line width = 1] (10,0)--(10,1);

\draw [line width = 1] (12,2)--(12,1);
\draw [line width = 1] (11,2)--(11,1);
\draw [line width = 1] (10,2)--(10,1);

\node at (.5,.5) {$+$};
\node at (1.5,.5) {$+$};
\node at (2.5,.5) {$+$};
\node at (10.5,.5) {$+$};
\node at (11.5,.5) {$+$};
\node at (12.5,.5) {$0$};

\node at (5.5,.5) {$+$};
\node at (6.5,.5) {$+$};
\node at (7.5,.5) {$+$};

\node at (1.5,1.5) {$-$};
\node at (2.5,1.5) {$-$};

\node at (10.5,1.5) {$-$};
\node at (11.5,1.5) {$-$};
\node at (12.5,1.5) {$-$};

\node at (5.5,1.5) {$-$};
\node at (6.5,1.5) {$-$};
\node at (7.5,1.5) {$-$};

\node at (.5,1.5) {$0$};

\node at (3.5,.5) {$\cdot$};
\node at (4.5,.5) {$\cdot$};

\node at (8.5,.5) {$\cdot$};
\node at (9.5,.5) {$\cdot$};

\node at (3.5,1.5) {$\cdot$};
\node at (4.5,1.5) {$\cdot$};

\node at (8.5,1.5) {$\cdot$};
\node at (9.5,1.5) {$\cdot$};

\node at (-.5, -.5) {$0$};
\end{tikzpicture}
\caption{Step 2 (-)}
\end{subfigure}
~~~~~~~~~~~~
\begin{subfigure}{.21\textwidth}
\begin{tikzpicture}[scale=.33]
\draw [line width = 1,<->] (0,14)--(0,0)--(14.5,0);

\draw [line width = 2] (1,0)--(14,0)--(14,1)--(2,1)--(2,2)--(1,2)--(1,0);
\draw [line width = 1] (0,2)--(13,2)--(13,0);
\draw [line width = 1] (0,1)--(1,1);
\draw [line width = 1] (1,0)--(1,1);
\draw [line width = 1] (2,0)--(2,1);
\draw [line width = 1] (3,0)--(3,1);
\draw [line width = 1] (1,1)--(2,1);

\draw [line width = 1] (1,2)--(1,1);
\draw [line width = 1] (2,2)--(2,1);
\draw [line width = 1] (3,2)--(3,1);
\draw [line width = 1] (12,1)--(13,1);

\draw [line width = 1] (5,2)--(5,1);
\draw [line width = 1] (6,2)--(6,1);
\draw [line width = 1] (7,2)--(7,1);
\draw [line width = 1] (8,2)--(8,1);
\draw [line width = 1] (4,2)--(4,1);
\draw [line width = 1] (9,2)--(9,1);

\draw [line width = 1] (5,0)--(5,1);
\draw [line width = 1] (6,0)--(6,1);
\draw [line width = 1] (7,0)--(7,1);
\draw [line width = 1] (8,0)--(8,1);
\draw [line width = 1] (4,0)--(4,1);
\draw [line width = 1] (9,0)--(9,1);

\draw [line width = 1] (12,0)--(12,1);
\draw [line width = 1] (11,0)--(11,1);
\draw [line width = 1] (10,0)--(10,1);

\draw [line width = 1] (12,2)--(12,1);
\draw [line width = 1] (11,2)--(11,1);
\draw [line width = 1] (10,2)--(10,1);

\node at (.5,.5) {$+$};
\node at (1.5,.5) {$0$};
\node at (2.5,.5) {$0$};
\node at (10.5,.5) {$0$};
\node at (11.5,.5) {$0$};
\node at (12.5,.5) {$-$};
\node at (13.5,.5) {$-$};

\node at (5.5,.5) {$0$};
\node at (6.5,.5) {$0$};
\node at (7.5,.5) {$0$};

\node at (1.5,1.5) {\tiny{-2}};
\node at (2.5,1.5) {$-$};

\node at (10.5,1.5) {$-$};
\node at (11.5,1.5) {$-$};
\node at (12.5,1.5) {$-$};

\node at (5.5,1.5) {$-$};
\node at (6.5,1.5) {$-$};
\node at (7.5,1.5) {$-$};

\node at (.5,1.5) {$0$};

\node at (3.5,.5) {$\cdot$};
\node at (4.5,.5) {$\cdot$};

\node at (8.5,.5) {$\cdot$};
\node at (9.5,.5) {$\cdot$};

\node at (3.5,1.5) {$\cdot$};
\node at (4.5,1.5) {$\cdot$};

\node at (8.5,1.5) {$\cdot$};
\node at (9.5,1.5) {$\cdot$};

\node at (-.5, -.5) {$0$};
\end{tikzpicture}
\caption{Step 3 (-)}
\end{subfigure}
~~~~~~~~~~~~~
\begin{subfigure}{.21\textwidth}
\begin{tikzpicture}[scale=.33]
\draw [line width = 1,<->] (0,14)--(0,0)--(14.5,0);

\draw [line width = 2] (12,0)--(14,0)--(14,1)--(13,1)--(13,13)--(12,13)--(12,0);
\draw [line width = 1] (0,2)--(13,2)--(13,0);
\draw [line width = 1] (0,1)--(1,1);
\draw [line width = 1] (1,0)--(1,1);
\draw [line width = 1] (2,0)--(2,1);
\draw [line width = 1] (3,0)--(3,1);
\draw [line width = 1] (1,1)--(13,1);

\draw [line width = 1] (12,3)--(13,3);
\draw [line width = 1] (12,4)--(13,4);
\draw [line width = 1] (12,5)--(13,5);
\draw [line width = 1] (12,6)--(13,6);
\draw [line width = 1] (12,7)--(13,7);
\draw [line width = 1] (12,8)--(13,8);
\draw [line width = 1] (12,9)--(13,9);
\draw [line width = 1] (12,10)--(13,10);
\draw [line width = 1] (12,11)--(13,11);
\draw [line width = 1] (12,12)--(13,12);

\draw [line width = 1] (1,2)--(1,1);
\draw [line width = 1] (2,2)--(2,1);
\draw [line width = 1] (3,2)--(3,1);
\draw [line width = 1] (12,1)--(13,1);

\draw [line width = 1] (5,2)--(5,1);
\draw [line width = 1] (6,2)--(6,1);
\draw [line width = 1] (7,2)--(7,1);
\draw [line width = 1] (8,2)--(8,1);
\draw [line width = 1] (4,2)--(4,1);
\draw [line width = 1] (9,2)--(9,1);

\draw [line width = 1] (5,0)--(5,1);
\draw [line width = 1] (6,0)--(6,1);
\draw [line width = 1] (7,0)--(7,1);
\draw [line width = 1] (8,0)--(8,1);
\draw [line width = 1] (4,0)--(4,1);
\draw [line width = 1] (9,0)--(9,1);

\draw [line width = 1] (12,0)--(12,1);
\draw [line width = 1] (11,0)--(11,1);
\draw [line width = 1] (10,0)--(10,1);

\draw [line width = 1] (12,2)--(12,1);
\draw [line width = 1] (11,2)--(11,1);
\draw [line width = 1] (10,2)--(10,1);

\node at (.5,.5) {$+$};
\node at (1.5,.5) {$0$};
\node at (2.5,.5) {$0$};
\node at (10.5,.5) {$0$};
\node at (11.5,.5) {$0$};
\node at (12.5,.5) {$0$};
\node at (13.5,.5) {$0$};

\node at (12.5,2.5) {$+$};
\node at (12.5,3.5) {$\cdot$};
\node at (12.5,4.5) {$\cdot$};
\node at (12.5,5.5) {$+$};
\node at (12.5,6.5) {$+$};
\node at (12.5,7.5) {$+$};
\node at (12.5,8.5) {$\cdot$};
\node at (12.5,9.5) {$\cdot$};
\node at (12.5,10.5) {$+$};
\node at (12.5,11.5) {$+$};
\node at (12.5,12.5) {$+$};

\node at (5.5,.5) {$0$};
\node at (6.5,.5) {$0$};
\node at (7.5,.5) {$0$};

\node at (1.5,1.5) {\tiny{-2}};
\node at (2.5,1.5) {$-$};

\node at (10.5,1.5) {$-$};
\node at (11.5,1.5) {$-$};
\node at (12.5,1.5) {$0$};

\node at (5.5,1.5) {$-$};
\node at (6.5,1.5) {$-$};
\node at (7.5,1.5) {$-$};

\node at (.5,1.5) {$0$};

\node at (3.5,.5) {$\cdot$};
\node at (4.5,.5) {$\cdot$};

\node at (8.5,.5) {$\cdot$};
\node at (9.5,.5) {$\cdot$};

\node at (3.5,1.5) {$\cdot$};
\node at (4.5,1.5) {$\cdot$};

\node at (8.5,1.5) {$\cdot$};
\node at (9.5,1.5) {$\cdot$};

\node at (-.5, -.5) {$0$};
\end{tikzpicture}
\caption{Step 4 (+)}
\end{subfigure}
~~~~~~~~~~~~~
\begin{subfigure}{.21\textwidth}
\begin{tikzpicture}[scale=.33]
\draw [line width = 1,<->] (0,14)--(0,0)--(14.5,0);

\draw [line width = 2] (12,0)--(14,0)--(14,1)--(13,1)--(13,13)--(12,13)--(12,0);
\draw [line width = 1] (0,2)--(13,2)--(13,0);
\draw [line width = 1] (0,1)--(1,1);
\draw [line width = 1] (1,0)--(1,1);
\draw [line width = 1] (2,0)--(2,1);
\draw [line width = 1] (3,0)--(3,1);
\draw [line width = 1] (1,1)--(13,1);

\draw [line width = 1] (12,3)--(13,3);
\draw [line width = 1] (12,4)--(13,4);
\draw [line width = 1] (12,5)--(13,5);
\draw [line width = 1] (12,6)--(13,6);
\draw [line width = 1] (12,7)--(13,7);
\draw [line width = 1] (12,8)--(13,8);
\draw [line width = 1] (12,9)--(13,9);
\draw [line width = 1] (12,10)--(13,10);
\draw [line width = 1] (12,11)--(13,11);
\draw [line width = 1] (12,12)--(13,12);

\draw [line width = 1] (1,2)--(1,1);
\draw [line width = 1] (2,2)--(2,1);
\draw [line width = 1] (3,2)--(3,1);
\draw [line width = 1] (12,1)--(13,1);

\draw [line width = 1] (5,2)--(5,1);
\draw [line width = 1] (6,2)--(6,1);
\draw [line width = 1] (7,2)--(7,1);
\draw [line width = 1] (8,2)--(8,1);
\draw [line width = 1] (4,2)--(4,1);
\draw [line width = 1] (9,2)--(9,1);

\draw [line width = 1] (5,0)--(5,1);
\draw [line width = 1] (6,0)--(6,1);
\draw [line width = 1] (7,0)--(7,1);
\draw [line width = 1] (8,0)--(8,1);
\draw [line width = 1] (4,0)--(4,1);
\draw [line width = 1] (9,0)--(9,1);

\draw [line width = 1] (12,0)--(12,1);
\draw [line width = 1] (11,0)--(11,1);
\draw [line width = 1] (10,0)--(10,1);

\draw [line width = 1] (12,2)--(12,1);
\draw [line width = 1] (11,2)--(11,1);
\draw [line width = 1] (10,2)--(10,1);

\draw [line width = 1] (2,2)--(2,3)--(3,3)--(3,4)--(4,4)--(4,5)--(5,5)--(5,6)--(6,6)--(6,7)--(7,7)--(7,8)--(8,8)--(8,9)--(9,9)--(9,10)--(10,10)--(10,11)--(11,11)--(11,12)--(12,12)--(12,13);

\draw [line width = 1] (2,2)--(3,2)--(3,3)--(4,3)--(4,4)--(5,4)--(5,5)--(6,5)--(6,6)--(7,6)--(7,7)--(8,7)--(8,8)--(9,8)--(9,9)--(10,9)--(10,10)--(11,10)--(11,11)--(12,11)--(12,12)--(13,12);

\node at (6,-.75) {$k-1$};
\draw [fill] (6,0) circle (4pt);

\node at (1.5,6) {$k-1$};
\draw [fill] (0,6) circle (4pt);

\node at (12,-.75) {$2k-2$};
\draw [fill] (12,0) circle (4pt);

\node at (1.8,12) {$2k-2$};
\draw [fill] (0,12) circle (4pt);

\node at (2.5,2.5) {$0$};
\node at (3.5,3.5) {$\cdot$};
\node at (4.5,4.5) {$\cdot$};
\node at (5.5,5.5) {$0$};
\node at (6.5,6.5) {$-$};
\node at (7.5,7.5) {$0$};
\node at (8.5,8.5) {$\cdot$};
\node at (9.5,9.5) {$\cdot$};
\node at (10.5,10.5) {$0$};
\node at (11.5,11.5) {$0$};

\node at (.5,.5) {$+$};
\node at (1.5,.5) {$0$};
\node at (2.5,.5) {$0$};
\node at (10.5,.5) {$0$};
\node at (11.5,.5) {$0$};
\node at (12.5,.5) {$0$};
\node at (13.5,.5) {$0$};

\node at (12.5,2.5) {$0$};
\node at (12.5,3.5) {$\cdot$};
\node at (12.5,4.5) {$\cdot$};
\node at (12.5,5.5) {$0$};
\node at (12.5,6.5) {$0$};
\node at (12.5,7.5) {$0$};
\node at (12.5,8.5) {$\cdot$};
\node at (12.5,9.5) {$\cdot$};
\node at (12.5,10.5) {$0$};
\node at (12.5,11.5) {$0$};
\node at (12.5,12.5) {$+$};

\node at (5.5,.5) {$0$};
\node at (6.5,.5) {$0$};
\node at (7.5,.5) {$0$};

\node at (1.5,1.5) {$-$};
\node at (2.5,1.5) {$0$};

\node at (10.5,1.5) {$0$};
\node at (11.5,1.5) {$0$};
\node at (12.5,1.5) {$0$};

\node at (5.5,1.5) {$0$};
\node at (6.5,1.5) {$0$};
\node at (7.5,1.5) {$0$};

\node at (.5,1.5) {$0$};

\node at (3.5,.5) {$\cdot$};
\node at (4.5,.5) {$\cdot$};

\node at (8.5,.5) {$\cdot$};
\node at (9.5,.5) {$\cdot$};

\node at (3.5,1.5) {$\cdot$};
\node at (4.5,1.5) {$\cdot$};

\node at (8.5,1.5) {$\cdot$};
\node at (9.5,1.5) {$\cdot$};

\node at (-.5, -.5) {$0$};
\end{tikzpicture}
\caption{Step 5}
\end{subfigure}
\caption{Building $C_5(k)$ out of $\{H_1(k), H_2(k), H_3(k), H_4(k)\}.$}
\label{fig:main-diag345}
\end{figure}

\begin{figure}[h!]
\centering
\begin{tikzpicture}[scale=.4]
\draw [line width = 1,<->] (0,5)--(0,0)--(5,0);
\draw [line width = 1] (0,0)--(1,0)--(1,1)--(2,1)--(2,2)--(3,2)--(3,3)--(4,3)--(4,4)--(3,4)--(3,3)--(2,3)--(2,2)--(1,2)--(1,1)--(0,1)--(0,0);

\draw [line width = 1,<->] (7,5)--(7,0)--(12,0);
\draw [line width = 1] (7,0)--(8,0)--(8,1)--(9,1)--(9,2)--(10,2)--(10,3)--(11,3)--(11,4)--(10,4)--(10,3)--(9,3)--(9,2)--(8,2)--(8,1)--(7,1)--(7,0);
\draw [line width = 1] (9,1)--(10,1)--(10,2)--(11,2)--(11,3);

\draw [line width = 1,<->] (14,5)--(14,0)--(19,0);
\draw [line width = 1] (14,0)--(15,0)--(15,1)--(16,1)--(16,2)--(17,2)--(17,3)--(17,3)--(18,3)--(18,3)--(17,3)--(16,3)--(16,2)--(15,2)--(15,1)--(14,1)--(14,0);
\draw [line width = 1] (16,1)--(17,1)--(17,2)--(18,2)--(18,3);

\draw [line width = 1,<->] (21,5)--(21,0)--(26,0);
\draw [line width = 1] (21,0)--(22,0)--(22,1)--(23,1)--(23,2)--(24,2)--(24,3)--(23,3)--(23,2)--(22,2)--(22,1)--(21,1)--(21,0);

\node at (15.5,1.5) {$-$};
\node at (16.5,1.5) {$-$};
\node at (16.5,2.5) {$+$};
\node at (17.5,2.5) {$+$};

\node at (9.5,2.5) {$+$};
\node at (10.5,3.5) {$-$};
\node at (9.5,1.5) {$+$};
\node at (10.5,2.5) {$-$};
\node at (6.5, -.5) {$0$};
\node at (13.5, -.5) {$0$};
\node at (20.5, -.5) {$0$};

\node at (23.5, 2.5) {\tiny{+2}};
\node at (22.5, 1.5) {$-$};

\node at (5.5,2.5) {$+$};
\node at (12.5,2.5) {$+$};
\node at (19.5,2.5) {$=$};

\node at (3.5,3.5) {$+$};
\node at (-.5, -.5) {$0$};
\end{tikzpicture}
\caption{Tiles arithmetic: $x^3+y^3-x^2yC_4(k)+xyC_3(k)=2x^2y^2-xy$.}
\label{fig:diagonaltiles-bis22}
\end{figure}
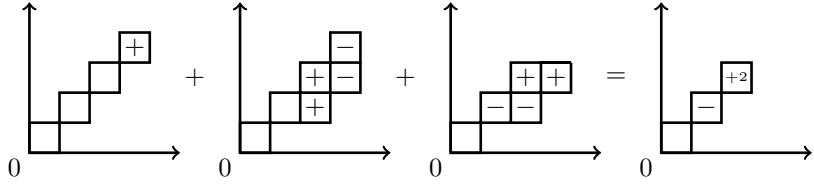

\end{pf}

\begin{propo} The sets $\{C_i(k), 1\le i\le 5\}$ and $\{H_i(k),1\le i\le 4\}$ generate the same ideal in $\mathbb{Z}[X,Y]$.
\end{propo}

\begin{pf} This follows from Propositions~\ref{p:prop5-even},~\ref{p:prop6-even}.
\end{pf}

\begin{propo} One has the following formulas:
\begin{equation}
\begin{aligned}
&S(C_1(k),C_2(k))=-y^{k-1}C_1(k)+x^{k-1}C_2(k)-y^{k-1}(1+x^2+\cdots +x^{k-3})C_3(k)\\
&+x^{k-1}(1+y^2+\cdots +y^{k-3})C_4(k)+y^{k-1}C_5(k)-y^{k-1}\left \lfloor{\frac{k-1}{2}}\right \rfloor C_4(k)-x^{k-1}C_5(k)+x^{k-1}\left \lfloor{\frac{k-1}{2}}\right \rfloor C_3(k),
\end{aligned}
\end{equation}

\begin{equation}\label{eq:c1-c3}
\begin{aligned}
S(C_1(k),C_3(k))&=xC_2(k)-y^{k-2}C_4(k)+y^{k-2}C_3(k)+(xy^{k-4}+xy^{k-6}+\dots +xy+\epsilon(k))C_4(k)\\
&-xC_5(k)+\left \lfloor{\frac{k-1}{2}}\right \rfloor xC_3(k),
\end{aligned}
\end{equation}

\begin{equation}\label{eq:c1-c4}
\begin{aligned}
S(C_1(k),C_4(k))&=C_2(k)+(x^{k-4}+x^{k-6}+\dots +x^2+\epsilon(k))C_4(k)-C_5(k)+\left \lfloor{\frac{k-1}{2}}\right \rfloor C_3(k),
\end{aligned}
\end{equation}

\begin{equation}\label{eq:c1-c5}
\begin{aligned}
S(C_1(k),C_5(k))&=(k-2)C_2(k)+(k-2)C_3(k)(1-\epsilon(k)+y+y^3+\dots +y^{k-3})\\
&+\left (2\left \lfloor{\frac{k-3}{2}}\right \rfloor+\epsilon(k)\right ) C_5(k)
+(k-2)\left \lfloor{\frac{k-1}{2}}\right \rfloor C_3(k),
\end{aligned}
\end{equation}

\begin{equation}\label{eq:c2-c3}
\begin{aligned}
S(C_2(k),C_3(k))&=C_1(k)+(y^{k-4}+y^{k-6}+\dots +y^2+\epsilon(k))C_3(k)-C_5(k)+\left \lfloor{\frac{k-1}{2}}\right \rfloor C_4(k),
\end{aligned}
\end{equation}

\begin{equation}\label{eq:c2-c4}
\begin{aligned}
S(C_2(k),C_4(k))&=yC_1(k)-x^{k-2}C_3(k)+x^{k-2}C_4(k)+(x^{k-4}y+x^{k-6}y+\dots +x^2y+\epsilon(k)y)C_3(k)\\
&-yC_5(k)+\left \lfloor{\frac{k-1}{2}}\right \rfloor yC_4(k),
\end{aligned}
\end{equation}

\begin{equation}\label{eq:c2-c5}
\begin{aligned}
S(C_2(k),C_5(k))&=(k-2)C_1(k)+(k-2)C_4(k)(1-\epsilon(k)+x+x^3+\dots +x^{k-3})\\
&+\left (2\left \lfloor{\frac{k-3}{2}}\right \rfloor+\epsilon(k)\right ) C_5(k)+(k-2)\left \lfloor{\frac{k-1}{2}}\right \rfloor C_4(k),
\end{aligned}
\end{equation}

\begin{equation}\label{eq:c3-c4}
\begin{aligned}
S(C_3(k),C_4(k))&=-C_3(k)+C_4(k)
\end{aligned}
\end{equation}

\begin{equation}\label{eq:c3-c5}
\begin{aligned}
S(C_3(k),C_5(k))&=C_5(k),\ \ S(C_4(k),C_5(k))=C_5(k),
\end{aligned}
\end{equation}
where $\epsilon(k)=\frac{1-(-1)^k}{2}$, which are given by $D$-reductions. Therefore, $\{C_i(k),i\le 1\le 5\}$ form a Gr\"obner basis.
\end{propo}

\begin{pf} We observe that we can always choose one of the coefficients $c_1, c_2$ in Definition 3 to be zero. So in order to check that we have a Gr\"obner basis, we do not need to use $G$-polynomials.

Due to the symmetry, some formulas above follows immediately from others: \eqref{eq:c2-c3} follows from \eqref{eq:c1-c4}, \eqref{eq:c2-c4} follows from \eqref{eq:c1-c3}, \eqref{eq:c2-c5} follows from \eqref{eq:c1-c5}, and second formula in~\eqref{eq:c3-c5} follows from the first. For the rest, note that the leading monomial in $C_1(k)$ is $y^k$, the leading monomial in $C_2(k)$ is $x^k$, the leading monomial in $C_3(k)$ is $x^2y$, the leading monomial in $C_4(k)$ is $xy^2$, and the leading monomial in $C_5(k)$ is $(k-2)xy$.

The $D$-reduction of $S(C_1(k),C_2(k))$ is shown in Figure~\ref{fig:Dredc1c2}. $S(C_1(k),C_2(k))$ consists of two disjoint symmetric tiles. The reduction of them is similar and it is shown in parallel in Figure~\ref{fig:Dredc1c2}.
We start with
\begin{equation*}
S(C_1(k),C_2(k))=x^kC_1(k)-y^{k}C_2(k).
\end{equation*}

\begin{figure}[h!]
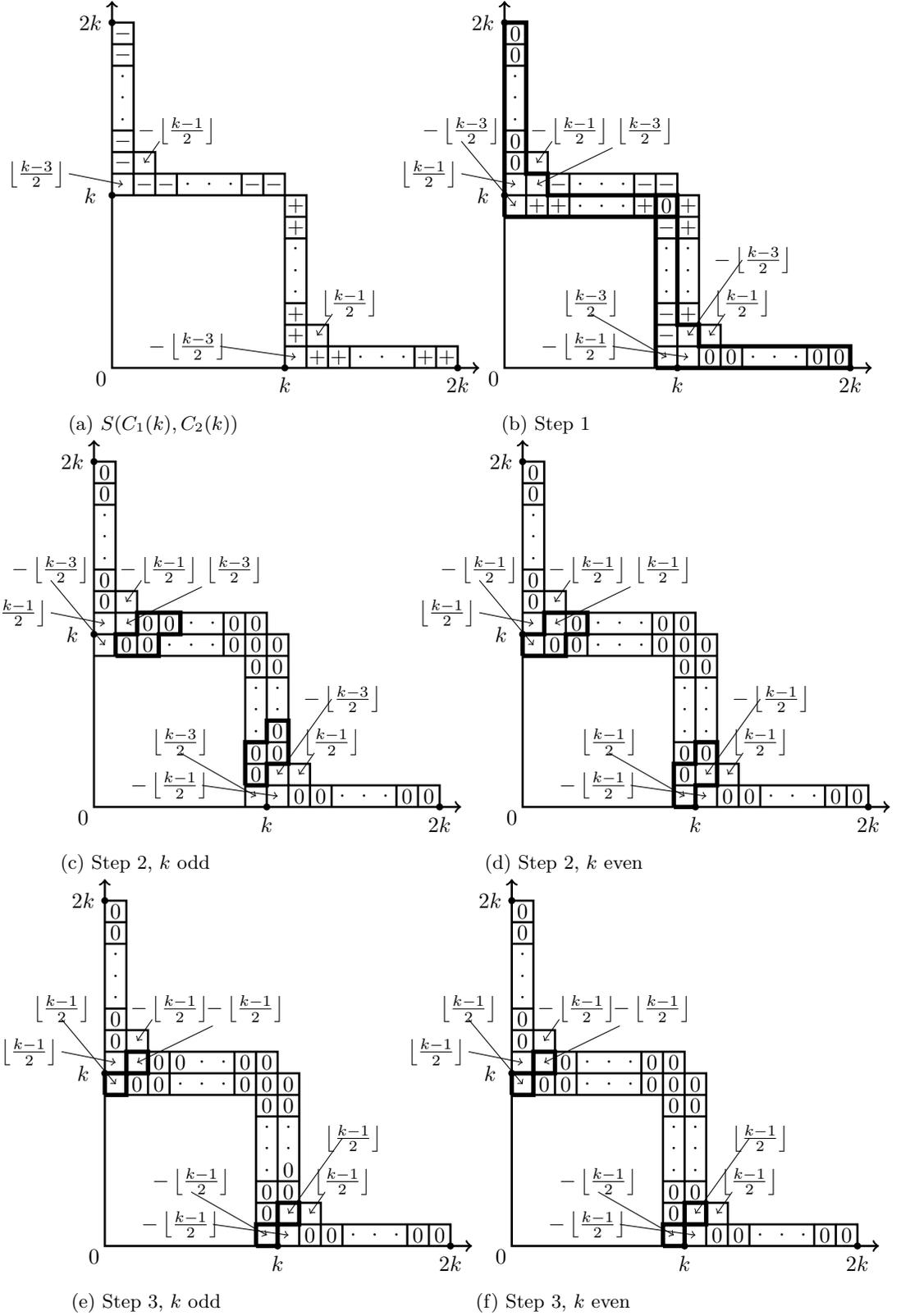

\centering
\begin{subfigure}{.3\textwidth}

\caption{$S(C_1(k),C_2(k))$}
\end{subfigure}
~~~~~~~~~~
\begin{subfigure}{.3\textwidth}
%
\caption{Step 1}
\end{subfigure}
~~~~~~~~~~~~
\begin{subfigure}{.3\textwidth}
%
\caption{Step 2, $k$ odd}
\end{subfigure}
~~~~~~~~~~~~~~~
\begin{subfigure}{.3\textwidth}
%
\caption{Step 2, $k$ even}
\end{subfigure}
~~~~~~~~~~~~~
\begin{subfigure}{.3\textwidth}
%
\caption{Step 3, $k$ odd}
\end{subfigure}
~~~~~~~~~~~~
\begin{subfigure}{.3\textwidth}
%
\caption{Step 3, $k$ even}
\end{subfigure}
\caption{The $D$-reduction of $S(C_1(k),C_2(k)).$}
\label{fig:Dredc1c2}
\end{figure}

The $D$-reduction of $S(C_1(k),C_3(k))$ is shown in Figure~\ref{fig:Dredc1c3}. We start with
\begin{equation*}
S(C_1(k),C_3(k))=x^2C_1(k)-y^{k-1}C_3(k).
\end{equation*}

From Step 3 to Step 4 we subtract $(xy^{k-4}+xy^{k-6}+\dots +xy)C_4(k)$ or $(xy^{k-4}+xy^{k-6}+\dots +xy^2+x)C_4(k)$, depending on $k$ odd or even. From Step 4 to Step 5 we use the following formulas:
\begin{equation}\label{eq:123odd}
(k-2)-\left \lfloor{\frac{k-3}{2}}\right \rfloor=\left \lfloor{\frac{k-1}{2}}\right \rfloor, \text{ if $k$ is odd },\ \
(k-2)-\left \lfloor{\frac{k-1}{2}}\right \rfloor=\left \lfloor{\frac{k-1}{2}}\right \rfloor, \text{ if $k$ is even }.
\end{equation}

The $D$-reduction of $S(C_1(k),C_4(k))$ is shown in Figure~\ref{fig:Dredc1c4}. We start with
\begin{equation}\label{eq:123even}
S(C_1(k),C_4(k))=xC_1(k)-y^{k-2}C_4(k).
\end{equation}

From Step 1 to Step 2 we subtract $(x^{k-4}+x^{k-6}+\dots +x)C_4(k)$ or $(x^{k-4}+x^{k-6}+\dots +x^2+1)C_4(k)$, depending on $k$ odd or even. From Step 2 to Step 3 we use formulas~\ref{eq:123odd}.

\begin{figure}[h!]
\centering
\begin{subfigure}{.35\textwidth}

\caption{$S(C_1(k),C_3(k))$}
\end{subfigure}
~~~~~~~~~~
\begin{subfigure}{.35\textwidth}
%
\caption{Step 1 (-)}
\end{subfigure}
~~~~~~~~~~~~
\begin{subfigure}{.35\textwidth}
%
\caption{Step 2 (+)}
\end{subfigure}
~~~~~~~~~~~~
\begin{subfigure}{.35\textwidth}
%
\caption{Step 3 (-)}
\end{subfigure}
~~~~~~~~~~~~
\begin{subfigure}{.35\textwidth}
%
\caption{Step 4 (-), $k$ odd}
\end{subfigure}
~~~~~~~~~~~
\begin{subfigure}{.35\textwidth}
%
\caption{Step 4 (-), $k$ even}
\end{subfigure}
~~~~~~~~~~~
\begin{subfigure}{.35\textwidth}
%
\caption{Step 5 (+), $k$ odd}
\end{subfigure}
~~~~~~~~~~~
\begin{subfigure}{.35\textwidth}
%
\caption{Step 5 (+), $k$ even}
\end{subfigure}
\caption{The $D$-reduction of $S(C_1(k),C_3(k)).$}
\label{fig:Dredc1c3}
\end{figure}

\begin{figure}[h!]
\centering
\begin{subfigure}{.35\textwidth}
%
\caption{$S(C_1(k),C_4(k))$}
\end{subfigure}
~~~~~~~~~~
\begin{subfigure}{.35\textwidth}
%
\caption{Step 1 (-)}
\end{subfigure}
~~~~~~~~~~~~
\begin{subfigure}{.35\textwidth}
%
\caption{Step 2 (-), $k$ odd}
\end{subfigure}
~~~~~~~~~~~~
\begin{subfigure}{.35\textwidth}
%
\caption{Step 2 (-), $k$ even}
\end{subfigure}
~~~~~~~~~~~~
\begin{subfigure}{.35\textwidth}
%
\caption{Step 3 (-), $k$ odd}
\end{subfigure}
~~~~~~~~~~~~
\begin{subfigure}{.35\textwidth}
%
\caption{Step 3 (-), $k$ even}
\end{subfigure}
\caption{The $D$-reduction of $S(C_1(k),C_4(k)).$}
\label{fig:Dredc1c4}
\end{figure}

\begin{figure}[h!]
\centering
\begin{subfigure}{.35\textwidth}
%
\caption{$S(C_1(k),C_5(k))$}
\end{subfigure}
~~~~~~~~~~
\begin{subfigure}{.35\textwidth}
%
\caption{Step 1, (-)}
\end{subfigure}
~~~~~~~~~~~~
\begin{subfigure}{.35\textwidth}
%
\caption{Step 2, (-), $k$ odd}
\end{subfigure}
~~~~~~~~~~~~
\begin{subfigure}{.35\textwidth}
%
\caption{Step 2, (-), $k$ even}
\end{subfigure}
~~~~~~~~~~~~
\begin{subfigure}{.35\textwidth}
%
\caption{Step 3, (+), $k$ odd}
\end{subfigure}
~~~~~~~~~~~~
\begin{subfigure}{.35\textwidth}
%
\caption{Step 3, (+), $k$ even}
\end{subfigure}
\caption{The $D$-reduction of $S(C_1(k),C_5(k)).$}
\label{fig:Dredc1c5}
\end{figure}

The $D$-reduction of $S(C_1(k),C_5(k))$ is shown in Figure~\ref{fig:Dredc1c5}. We start with
\begin{equation}\label{eq:123even}
S(C_1(k),C_5(k))=(k-2)xC_1(k)-y^{k-1}C_5(k).
\end{equation}

To reach Step 1, we subtract $(k-2)C_2(k)$. To reach Step 2, we subtract $(1+y^2+\dots +y^{k-3})C_4(k)$ if $k$ is odd and $(y+\dots +y^{k-3})C_4(k)$ if $k$ is
even. To reach Step 3, we add $2\left \lfloor{\frac{k-1}{2}}\right \rfloor C_5(k)$ if $k$ is odd and $\left ( 2\left \lfloor{\frac{k-3}{2}}\right \rfloor+1\right ) C_5(k)$ if $k$ is even.

The $D$-reduction of $S(C_3(k),C_4(k))$ is:
\begin{equation*}
\begin{gathered}
S(C_3(k),C_4(k))=yC_3(k)-xC_4(k)\\
=x^2y^2+xy^2-xy-y-(x^2y^2+x^2y-xy-x)=xy^2-y-x^2y+x=-C_3(k)+C_4(k).
\end{gathered}
\end{equation*}

The $D$-reduction of $S(C_3(k),C_5(k))$ is:

\begin{equation*}
\begin{gathered}
S(C_3(k),C_5(k))=(k-2)C_3(k)-xC_5(k)\\
=(k-2)x^2y+(k-2)xy-(k-2)x-(k-2)-(k-2)x^2y+(k-2)x=C_5(k).
\end{gathered}
\end{equation*}
\end{pf}

\section{Gr\"obner basis for $\mathcal{T}_n^+, n$ even}\label{s:3-bis}
We consider first the case $n=4$.

\begin{propo}\label{p:propt+} A Gr\"obner basis for the ideal generated by $\mathcal{T}_4^+$ is given by: $D_1=y^2 + 2y + 1, D_2=x - y$.
\end{propo}

\begin{pf} One has:
\begin{equation*}
\begin{gathered}
y^2+2y+1= (y^2 + y + 1 + x) + (x^2 + y + xy + x^2y)-x(xy+x+y+1)\\
x-y=-(x^2 + y + xy + x^2y)+x(xy+x+y+1),
\end{gathered}
\end{equation*}
thus the Gr\"obner basis can be generated by $\mathcal{T}_4^+$. Conversely, one has:
\begin{equation*}
\begin{gathered}
y^2 + y + 1 + x=(y^2 + 2y + 1)+(x-y),\\
y^2 +xy^2 + xy + x=y(y^2+y+1)+(y^2 + y + 1)(x-y),\\
y + 1 + x + x^2 =(y^2+2y+1)+(x + y + 1)(x-y),\\
x^2 + y + xy + x^2y=y(y^2+2y+1)+(xy + y^2 + x + 2y)(x-y),\\
xy+x+y+1=(y^2+2y+1)+(y + 1)(x-y),
\end{gathered}
\end{equation*}
thus $\mathcal{T}_4^+$ is generated by the Gr\"obner basis.

The $S$-polynomial $S(D_1,D_2)$ is reduced as follows:
\begin{equation*}
S(D_1,D_2)=x(y^2+2y+1)-y^2(x-y)=2xy+x+y^3= y(y^2+2y+1)+(2y + 1)(x-y).
\end{equation*}
\end{pf}

Let now $n=2k, k\ge 3$. Recall that $\mathcal{T}_n^{+}$ is the set $\mathcal{T}_n$ plus a $2\times 2$ tile with polynomial $H_5(k)=xy+x+y+1$. We show that the Gr\"obner basis for the ideal generated by $\mathcal{T}_n^{+}$ is given by:
\begin{equation*}
\begin{aligned}
D_1(k)&=y^2+2y+1,\\
D_2(k)&=\left[(k-2)(k-1) - 1\right]y + x + (k-2)(k-1),\\
D_3(n)&=k(k-2)\left(y+1\right).
\end{aligned}
\end{equation*}

\begin{figure}[h!]
\centering
\begin{subfigure}{.25\textwidth}
\begin{tikzpicture}[scale=.4]
\draw [line width = 1,<->] (0,4)--(0,0)--(2.5,0);
\draw [line width = 1] (0,3)--(1,3)--(1,0);
\draw [line width = 1] (0,2)--(1,2);
\draw [line width = 1] (0,1)--(1,1);

\node at (.5,.5) {$+$};
\node at (.5,1.5) {\tiny{$+2$}};
\node at (.5,2.5) {$+$};
\node at (-.5, -.5) {$0$};
\end{tikzpicture}
\caption{$D_1(k)$}
\end{subfigure}
~~~~~
\begin{subfigure}{.25\textwidth}
\begin{tikzpicture}[scale=.4]
\draw [line width = 1,<->] (0,3)--(0,0)--(3,0);
\draw [line width = 1] (0,2)--(1,2)--(1,1)--(2,1)--(2,0);
\draw [line width = 1] (0,1)--(1,1)--(1,0);

\node at (.5,.5) {};
\node at (.5,1.5) {};
\node at (1.5,.5) {$+$};
\node at (-.5, -.5) {$0$};

\node at (-2.7,1.5) {$(k-2)(k-1)$};
\draw [->] (-1.3,.9)--(.5,.5);

\node at (3.5,3.5) {$(k-2)(k-1)-1$};
\draw [->] (2,3)--(.5,1.5);
\end{tikzpicture}
\caption{$D_2(k)$}
\end{subfigure}
~~~~~
\begin{subfigure}{.25\textwidth}
\begin{tikzpicture}[scale=.4]
\draw [line width = 1,<->] (0,3)--(0,0)--(2.5,0);
\draw [line width = 1] (0,1)--(1,1);
\draw [line width = 1] (0,2)--(1,2)--(1,0);

\node at (.5,.5) {};
\node at (.5,1.5) {};
\node at (-.5, -.5) {$0$};

\node at (-2.5,1) {$k(k-2)$};
\draw [->] (-1,.8)--(.5,.5);

\node at (3,3) {$k(k-2)$};
\draw [->] (2,2.7)--(.5,1.5);
\end{tikzpicture}
\caption{$D_3(k)$}
\end{subfigure}
\caption{The Gr\"obner basis $\{D_1(k), D_2(k), D_3(k)\}.$}
\label{fig:tnpbasis}
\end{figure}
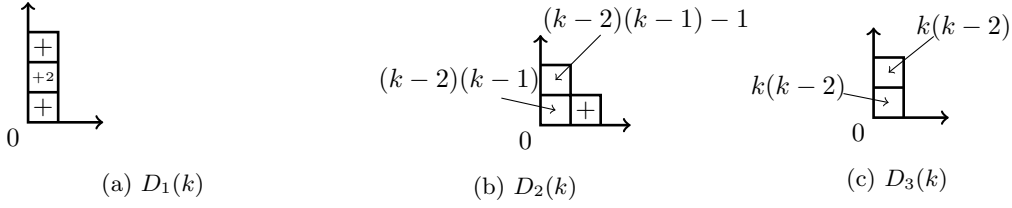

As $\mathcal{T}_n\subseteq \mathcal{T}_n^{+}$, $\mathcal{T}_n^+$ generates the Gr\"obner basis $C_i(k)$ for $\mathcal{T}_n$. Next formula shows how to generate $D_1(k)$:
\begin{equation}\label{eq:genP1}
\begin{aligned}
H_5(k)y - C_4(k)= (y+xy+y^2+xy^2) - (-1-y+xy+xy^2)= 1 + 2y + y^2=D_1(k).
\end{aligned}
\end{equation}

\begin{lem}\label{l:tnplem} The polynomial $1+x+(k-1)(y+y^2)$ is generated by $H_1(k)$ and $D_1(k),$ and $H_1(k)$ is generated by $1+x+(k-1)(y+y^2)$ and $D_1(k)$.
\end{lem}
\begin{pf} First produce $1+y-y^2-y^3 = D_1(k) (1-y)$. Adding copies of this tile to $H_1(k)$ gives the sum:
\begin{equation}\label{eq:some-pol-3030}
\begin{aligned}
H_1(k) + &\left [y^{n-5} + 2y^{n-7} + 3y^{n-9} + \ldots + (k-2)y\right ]\left(1+y-y^2-y^3\right ).
\end{aligned}
\end{equation}
Expanding \eqref{eq:some-pol-3030} gives a telescopic sum that reduces to $\left(y^2+y\right)\left(k-1\right) + 1 + x$.
\end{pf}

\begin{propo}\label{p:tnptogroebner} The polynomials $D_1(k)$, $D_2(k)$, and $D_3(k)$ belong to the ideal generated by $\mathcal{T}_n^{+}$.
\end{propo}
\begin{pf}
We showed in (\ref{eq:genP1}) how to generate $D_1(k)$. By Lemma~\ref{l:tnplem}, we can start from $H_1(k)$ to produce $1+x+\left(k-1\right)(y+y^2)$. Then, subtract as follows:
\begin{equation}\label{eq:some-pol0020202}
\begin{aligned}
1+x+\left(k-1\right)(y+y^2) - D_1(k)\left(k-1\right) &= -(k-2) + x - y\left(k-1\right).
\end{aligned}
\end{equation}

By the symmetry of $\mathcal{T}_n^{+}$ about $x=y$, we can also generate $-(k-2) - x\left(k-1\right) + y$. Combining the two:
\begin{equation*}
\begin{aligned}
\left(k-1\right)\left[-(k-2) + x - y\left(k-1\right)\right]+ \left(-(k-2) - x\left(k-1\right) + y\right) &= -D_3(k).
\end{aligned}
\end{equation*}

Finally, $D_2(k)$ is produced from $D_3(k)$ and $-(k-2) + x - y\left(k-1\right)$, which we have from~\eqref{eq:some-pol0020202}:
\begin{equation*}
\begin{aligned}
\left[-(k-2) + x - y\left(k-1\right)\right] + \left(y+1\right)k(k-2) &= D_2(k).
\end{aligned}
\end{equation*}

\end{pf}

\begin{lem}\label{l:tnplem2} The polynomials $\bar D_1(k)=1+2x+x^2$, $\bar D_2(k)=(k-1)(k-2)+x\left [(k-1)(k-2)-1\right ]+y$, and $\bar D_3(k)=\left(1+x\right)k(k-2)$ belong to the ideal generated by $D_1(k)$, $D_2(k)$, and $D_3(k)$.
\end{lem}
\begin{pf} We show independently in Proposition~\ref{p:groebnertotnp} below that $H_5(k)$ is also in this ideal. Then one can easily check that:
\begin{equation*}
\begin{gathered}
\bar D_1(k)=(1+x)D_2(k)-\left [(k-1)(k-2)-1\right ]H_5(k)\\
\bar D_2(k)=\left [(k-1)(k-2)-1\right ]D_2(k) - (k-1)(k-3)D_3(k)\\
\bar D_3(k)=k(k-2)D_2(k)-\left [(k-2)(k-1)-1\right ]D_3(k).
\end{gathered}
\end{equation*}
\end{pf}

\begin{propo}\label{p:groebnertotnp} The members of $\mathcal{T}_n^{+}$ belong to the ideal generated by $D_1(k)$, $D_2(k)$, and $D_3(k)$.
\end{propo}
\begin{pf}
One has after calculations:
\begin{equation*}
H_5(k)=(1+y)D_2(k)-\left [(k-2)(k-1)-1\right ]D_1(k).
\end{equation*}

To obtain $H_1(k)$, begin with
\begin{equation}\label{eq:t1fromgb}
D_2(k)-D_3(k)+(k-1) D_1(k)=1 + x + (k-1)(y+y^2).
\end{equation}

By Lemma~\ref{l:tnplem}, this tile may be transformed into $H_1(k)$ using only $D_1(k)$, $D_2(k)$, and $D_3(k)$. We also get $H_3(k)$ by symmetry in the following way. Swap the variables $x,y$ in Lemma~\ref{l:tnplem}. Then we have that $H_3(k)$ and $1+y+\left(k-1\right)\left(x+x^2\right)$ can each be produced from the other using either $\mathcal{T}_n^{+}$, which is symmetric about $x=y$, or the tiles $\bar D_1(k)$, $\bar D_2(k)$, and $\bar D_3(k)$. Then swapping $x, y$ in (\ref{eq:t1fromgb}) allows to obtain $1+y+\left(k-1\right)\left(x+x^2\right)$ from the `inverse basis' $\bar D_i(n)$, which in turn can be obtained from the Gr\"obner basis itself by Lemma~\ref{l:tnplem2}. Therefore the Gr\"obner basis also generates $H_3(k)$.

$C_4(k)=-1-y+xy+xy^2$ can be used to change $H_1(k)$ into $H_2(k)$ as follows:
\begin{equation*}
H_1(k)+(-1-y+xy+xy^2)(1+y^2+y^4+\ldots+y^{n-4}) = H_2(k).
\end{equation*}

By symmetry, the same process will change $H_3(k)$ into $H_4(k)$ using $C_3(k)$. It remains, then, to show that the Gr\"obner basis for $\mathcal{T}_n^{+}$ can generate $C_3(k)$ and $C_4(k)$.

Start with
\begin{equation*}
D_1(k)+y\left(D_3(k)-D_2(k)\right)= (y^2+y)k+1-xy.
\end{equation*}

Then, multiply by $y+1$ and add $-kyD_1(k)$:
\begin{equation*}
(y+1)\left[(y^2+y)k+1-xy\right] - ky\cdot(y^2+2y+1) = 1 + y - xy - xy^2=-C_4(k).
\end{equation*}

Once again, symmetry immediately gives us a corresponding procedure for $C_3(k)=1+x-xy-x^2y$, and the proof is complete.
\end{pf}

\begin{propo} The sets $\{D_i(k),1\le i\le 3\}$ and $\{H_i(k),1\le i\le 5\}$ generate the same ideal in $\mathbb{Z}[X,Y]$.
\end{propo}

\begin{pf} This follows from Propositions~\ref{p:tnptogroebner},~\ref{p:groebnertotnp}.
\end{pf}

\begin{propo} We have the following formulas:
\begin{align*}
S(D_1(k),D_2(k)) &= k(k-2)D_2(k) + (1-x)D_3(k),\\
S(D_1(k),D_3(k)) &= D_3(k),\\
S(D_2(k),D_3(k)) &= k(k-2)\left[1-(k-1)(k-2)\right]D_1(k) + k(k-2)D_2(k) - D_3(k),
\end{align*}
which are given by $D$-reductions. Therefore, $\{D_i(k),1\le i\le 3\}$ forms a Gr\"obner basis.
\end{propo}

\begin{pf} We start with
\begin{align*}
S(D_1(k),D_2(k)) &= k(k-2)\left[(k-2)(k-1)-1\right]D_1(k) - k(k-2)yD_2(k)\\
S(D_1(k),D_3(k)) &= k(k-2)D_1(k) - (y+1)D_3(k)\\
S(D_2(k),D_3(k)) &= -k(k-2)yD_2(k)+xD_3(k).
\end{align*}
The reader may easily check that the given reductions are valid for these $S$-polynomials.
\end{pf}

\section{Proof of Theorem~\ref{thm:main}}\label{s:4}

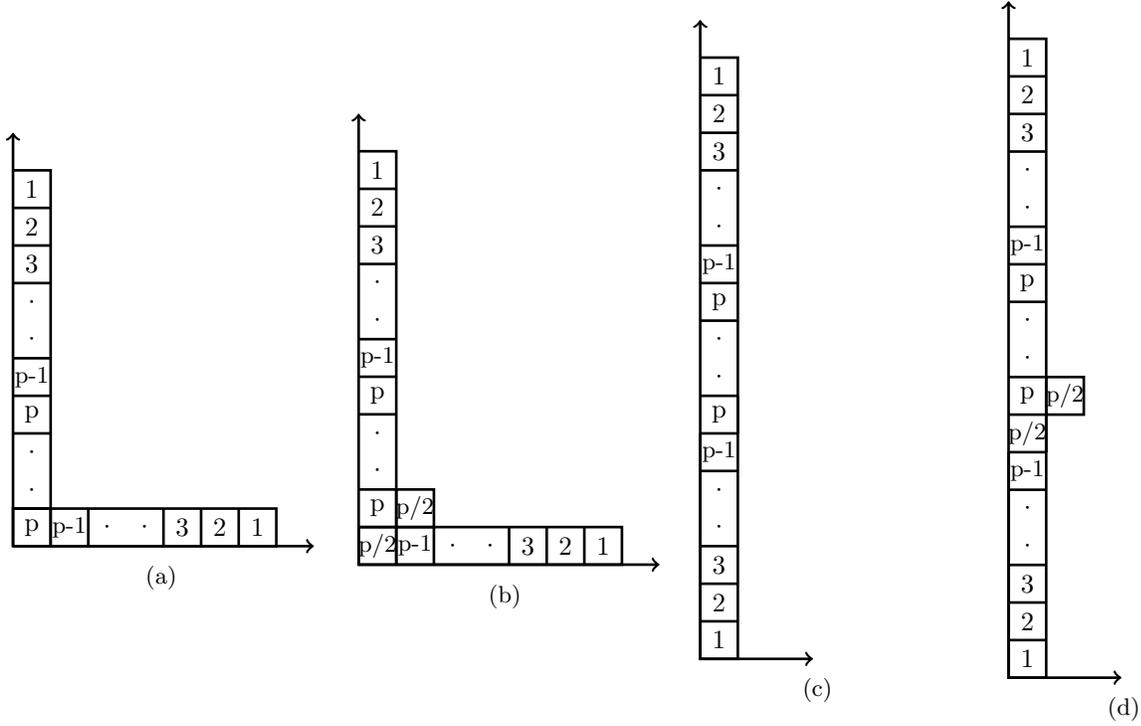
\begin{figure}[h]

\begin{subfigure}{.25\textwidth}
\begin{tikzpicture}[scale=.5]
\draw [line width = 1, <->](0,11)--(0,0)--(8,0);
\draw [line width = 1] (0,0)--(7,0)--(7,1)--(0,1)--(0,0);
\draw [line width = 1] (0,10)--(1,10)--(1,0);
\draw [line width = 1] (0,1)--(1,1);
\draw [line width = 1] (0,3)--(1,3);
\draw [line width = 1] (0,5)--(1,5);
\draw [line width = 1] (0,4)--(1,4);
\draw [line width = 1] (0,7)--(1,7);
\draw [line width = 1] (0,8)--(1,8);
\draw [line width = 1] (0,9)--(1,9);
\draw [line width = 1] (1,0)--(1,1);
\draw [line width = 1] (2,0)--(2,1);
\draw [line width = 1] (4,0)--(4,1);
\draw [line width = 1] (5,0)--(5,1);
\draw [line width = 1] (6,0)--(6,1);

\node at (.5, 9.5) {1};
\node at (.5, 8.5) {2};
\node at (.5, 7.5) {3};
\node at (.5, 6.5) {$\cdot$};
\node at (.5, 5.5) {$\cdot$};
\node at (.5, 4.5) {\small{p-1}};
\node at (.5, 3.5) {p};
\node at (.5, 2.5) {$\cdot$};
\node at (.5, 1.5) {$\cdot$};
\node at (.5, .5) {p};

\node at (6.5, .5) {1};
\node at (5.5, .5) {2};
\node at (4.5, .5) {3};
\node at (3.5, .5) {$\cdot$};
\node at (2.5, .5) {$\cdot$};
\node at (1.5, .5) {\small{p-1}};

\end{tikzpicture}
\caption{}
\end{subfigure}
~~
\begin{subfigure}{.25\textwidth}
\begin{tikzpicture}[scale=.5]
\draw [line width = 1, <->](0,12)--(0,0)--(8,0);
\draw [line width = 1] (0,0)--(7,0)--(7,1)--(0,1)--(0,0);
\draw [line width = 1] (0,11)--(1,11)--(1,0);
\draw [line width = 1] (0,1)--(1,1);

\draw [line width = 1] (0,2)--(1,2);
\draw [line width = 1] (0,4)--(1,4);
\draw [line width = 1] (0,6)--(1,6);
\draw [line width = 1] (0,5)--(1,5);
\draw [line width = 1] (0,8)--(1,8);
\draw [line width = 1] (0,9)--(1,9);
\draw [line width = 1] (0,10)--(1,10);

\draw [line width = 1] (1,0)--(1,1);
\draw [line width = 1] (2,0)--(2,1);
\draw [line width = 1] (4,0)--(4,1);
\draw [line width = 1] (5,0)--(5,1);
\draw [line width = 1] (6,0)--(6,1);

\draw [line width = 1] (1,1)--(2,1)--(2,2)--(1,2)--(1,1);

\node at (.5, 10.5) {1};
\node at (.5, 9.5) {2};
\node at (.5, 8.5) {3};
\node at (.5, 7.5) {$\cdot$};
\node at (.5, 6.5) {$\cdot$};
\node at (.5, 5.5) {\small{p-1}};
\node at (.5, 4.5) {p};
\node at (.5, 3.5) {$\cdot$};
\node at (.5, 2.5) {$\cdot$};
\node at (.5, 1.5) {p};

\node at (.5, .5) {\small{p/2}};
\node at (1.5, 1.5) {\small{p/2}};

\node at (6.5, .5) {1};
\node at (5.5, .5) {2};
\node at (4.5, .5) {3};
\node at (3.5, .5) {$\cdot$};
\node at (2.5, .5) {$\cdot$};
\node at (1.5, .5) {\small{p-1}};

\end{tikzpicture}
\caption{}
\end{subfigure}
~~
\begin{subfigure}{.2\textwidth}
\begin{tikzpicture}[scale=.5]
\draw [line width = 1, <->](0,11)--(0,-6)--(3,-6);
\draw [line width = 1] (0,0)--(0,-6)--(1,-6)--(1,0)--(0,0);
\draw [line width = 1] (0,10)--(1,10)--(1,0);
\draw [line width = 1] (0,1)--(1,1);
\draw [line width = 1] (0,3)--(1,3);
\draw [line width = 1] (0,5)--(1,5);
\draw [line width = 1] (0,4)--(1,4);
\draw [line width = 1] (0,7)--(1,7);
\draw [line width = 1] (0,8)--(1,8);
\draw [line width = 1] (0,9)--(1,9);
\draw [line width = 1] (1,0)--(1,1);
\draw [line width = 1] (0,-1)--(1,-1);
\draw [line width = 1] (0,-3)--(1,-3);
\draw [line width = 1] (0,-4)--(1,-4);
\draw [line width = 1] (0,-5)--(1,-5);

\node at (.5, 9.5) {1};
\node at (.5, 8.5) {2};
\node at (.5, 7.5) {3};
\node at (.5, 6.5) {$\cdot$};
\node at (.5, 5.5) {$\cdot$};
\node at (.5, 4.5) {\small{p-1}};
\node at (.5, 3.5) {p};
\node at (.5, 2.5) {$\cdot$};
\node at (.5, 1.5) {$\cdot$};
\node at (.5, .5) {p};

\node at (.5, -5.5) {1};
\node at (.5, -4.5) {2};
\node at (.5, -3.5) {3};
\node at (.5, -2.5) {$\cdot$};
\node at (.5, -1.5) {$\cdot$};
\node at (.5, -.5) {\small{p-1}};

\end{tikzpicture}
\caption{}
\end{subfigure}
~~~~~
\begin{subfigure}{.2\textwidth}
\begin{tikzpicture}[scale=.5]
\draw [line width = 1, <->](0,12)--(0,-6)--(3,-6);
\draw [line width = 1] (0,0)--(0,-6)--(1,-6)--(1,0)--(0,0);
\draw [line width = 1] (0,11)--(1,11)--(1,0);
\draw [line width = 1] (0,1)--(1,1);

\draw [line width = 1] (0,2)--(1,2);
\draw [line width = 1] (0,4)--(1,4);
\draw [line width = 1] (0,6)--(1,6);
\draw [line width = 1] (0,5)--(1,5);
\draw [line width = 1] (0,8)--(1,8);
\draw [line width = 1] (0,9)--(1,9);
\draw [line width = 1] (0,10)--(1,10);

\draw [line width = 1] (1,0)--(1,1);
\draw [line width = 1] (0,-1)--(1,-1);
\draw [line width = 1] (0,-3)--(1,-3);
\draw [line width = 1] (0,-4)--(1,-4);
\draw [line width = 1] (0,-5)--(1,-5);

\draw [line width = 1] (1,1)--(2,1)--(2,2)--(1,2)--(1,1);

\node at (.5, 10.5) {1};
\node at (.5, 9.5) {2};
\node at (.5, 8.5) {3};
\node at (.5, 7.5) {$\cdot$};
\node at (.5, 6.5) {$\cdot$};
\node at (.5, 5.5) {\small{p-1}};
\node at (.5, 4.5) {p};
\node at (.5, 3.5) {$\cdot$};
\node at (.5, 2.5) {$\cdot$};
\node at (.5, 1.5) {p};

\node at (.5, .5) {\small{p/2}};
\node at (1.5, 1.5) {\small{p/2}};

\node at (.5, -5.5) {1};
\node at (.5, -4.5) {2};
\node at (.5, -3.5) {3};
\node at (.5, -2.5) {$\cdot$};
\node at (.5, -1.5) {$\cdot$};
\node at (.5, -.5) {\small{p-1}};

\end{tikzpicture}
\caption{}
\end{subfigure}
\caption{$D$-reductions of a rectangle.}
\label{fig:newcong83}
\end{figure}

The case $n=6$ follows as in~\cite{nitica-L-shaped}. We assume for the rest of this section $n\ge 8$. Consider a $q\times p, q\ge p\ge 1,$ rectangle. Using the presence of $C_3(k)$ and $C_4(k)$ in the Gr\" obner basis, the rectangle can be reduced to one of the configurations in Figure~\ref{fig:newcong83}, a), b). Configuration b) appears when $q,p$ are both even. The number of cells labeled by $p$ is $q-p+1$ in a) and $q-p$ in b).

In what follows the signed tile $B=xy-1$ will play an important role. We recall that it can be moved horizontally/vertically as shown in Figure~\ref{fig:diagonaltiles}. The tile $B$ does not belong to the ideal generated by $\mathcal{T}_n$. Other signed tile of interest in the sequel is $D=y^{n+1}+y^n+y^{n-1}+\dots +y^2+y+1-xy$, which is the concatenation of a vertical bar of length $n$ and $B$. The tile $D=yH_1(k)-C_4(k)$ belongs to the ideal generated by $\mathcal{T}_n$.

Multiplying the polynomial associated to the rectangle by $y^{p}$, we can assume that the configurations in Figure~\ref{fig:newcong83} are at height $p-1$ above the $x$-axis. Using the tiles $C_3(k), C_4(k)$ and an  amount of tiles $B$ ($p/2$ if $p$ is even and zero if $p$ is odd), they can be reduced further to the configurations shown in Figure~\ref{fig:newcong83},c), d). We observe that b) is the sum of a) with $p/2$ copies of $B$.

Reducing further the configurations in Figure~\ref{fig:newcong83}, c), d), with copies of $D$, the existence of a signed tiling for the $q\times p$ rectangle becomes equivalent to deciding when the following two conditions are both true:

1) The the polynomial $Q(x)=1+y+y^2+\dots +y^{n-1}$ divides:
\begin{equation*}
P_{p,q}(y)=1+2y+3y^2+\dots +py^{p-1}+py^{p}+\dots +py^{q-1}+(p-1)y^q+(p-2)y^{q+1} +\dots +2y^{p+q-3}+y^{p+q-2}.
\end{equation*}

2) The extra tiles $B$ that appear while doing tile arithmetic for 1), including those from Figure~\ref{fig:newcong83}, can be cancelled out by $C_5(k)$.

If $p+q-1<n$, then $\deg Q>\deg P_{p,q}$, so divisibility does not hold. If $p+q-1\ge n$, we look at $P_{p,q}$ as a sum of $p$ polynomials with all coefficients equal to 1:
\begin{equation*}
\begin{aligned}
P_{p,q}(y)=1+y+y^2+y^3+\dots +y^{p-1}+y^{p}+\dots &+y^{q-1}+y^q+y^{q+1} +\dots +y^{p+q-4}+y^{p+q-3}+y^{p+q-2}\\
          +y+y^2+y^3+\dots +y^{p-1}+y^{p}+\dots &+y^{q-1}+y^q+y^{q+1} +y^{p+q-4}+\dots +y^{p+q-3}\\
          +y^3+\dots +y^{p-1}+y^{p}+\dots &+y^{q-1}+y^q+y^{q+1} +\dots +y^{p+q-4}\\
          \dots \dots &\dots \dots\\
        +y^{p}+\dots &+y^{q-1}.
\end{aligned}
\end{equation*}

We discuss first 1) and show that it is true when $p$ or $q$ is divisible by $n$. Then, assuming this condition satisfied, we discuss 2).

1) Assume that $p+q-1=nm+r, 0\le r<n,$ and $p=ns+t, 0\le t<n.$ The remainder $R_{p,q}(y)$ of the division of $P_{p,q}(y)$ by $Q(y)$ is the sum of the remainders of the division of the $p$ polynomials above by $Q(y)$.

If $r$ is odd, one has the sequence of remainders, each remainder written in a separate pair of parentheses:
\begin{equation*}
\begin{aligned}
R_{p,q}(y)=&(1+y+y^2+\dots +y^{r-1})+(y+y^2+\dots +y^{r-2})+(y^2+\dots +y^{r-3})+\cdots \cdots \\
+&(y^{\frac{r-1}{2}})-(y^{\frac{r-1}{2}})- \cdots \cdots \cdots \\
-&(y+y^2+\dots +y^{r-2})-(1+y+y^2+\dots +y^{r-1})\\
+&(y^{r+1}+y^{r+2}+\dots +y^{n-3}+y^{n-2})+(y^{r+2}+\dots +y^{n-3})+\cdots \cdots \cdots\\
+&(y^{\frac{r+n-1}{2}})-(y^{\frac{r+n-1}{2}})-\cdots \cdots \cdots\\
-&(y^{r+2}+\dots +y^{n-3})-(y^{r+1}+y^{r+3}+\dots +y^{n-3}+y^{n-2})+\cdots \cdots \cdots
\end{aligned}
\end{equation*}

If $p\ge n$, the sequence of remainders above is periodic with period $n$, given by the part of the sequence shown above, and the sum of any subsequence of $n$ consecutive remainders is 0. So if $p$ is divisible by $n$, $P_{p,q}(y)$ is divisible by $Q(y)$. If $p$ is not divisible by $n$, then doing first the cancellation as above and then using the symmetry present in the sequence of remainders, the sum of the sequence of remainders equals 0 only if $r+1=t$, that is, only if $q$ is divisible by $n$.

If $r$ is even, one has the sequence of remainders, each remainder written in a separate pair of parentheses:
\begin{equation*}
\begin{aligned}
R_{p,q}(y)=&(1+y+y^2+\dots +y^{r-1})+(y+y^2+\dots +y^{r-2})+(y^2+\dots +y^{r-3})+\cdots \cdots\\
+&(y^{\frac{r-2}{2}}+y^{\frac{r}{2}})+(0)-(y^{\frac{r-2}{2}}+y^{\frac{r}{2}})-\cdots\cdots\\
&-(y+y^2+\dots +y^{r-2})-(1+y+y^2+\dots +y^{r-1})\\
+&(y^{r+1}+y^{r+3}+\dots +y^{n-3}+y^{n-2})+(y^{r+2}+\dots +y^{n-3})+\cdots \cdots \cdots\\
+&(y^{\frac{r+n-1}{2}}+y^{\frac{r+n+1}{2}})+(0)-(y^{\frac{r+n-1}{2}}+y^{\frac{r+n+1}{2}})-\cdots \cdots \cdots\\
&-(y^{r+2}+\dots +y^{n-3})-(y^{r+1}+y^{r+3}+\dots +y^{n-3}+y^{n-2})-\cdots \cdots \cdots
\end{aligned}
\end{equation*}

If $p\ge n$, the sequence of remainders above is periodic with period $n$, given by the part of the sequence shown above, and the sum of any subsequence of $n$ consecutive remainders is 0. So if $p$ is divisible by $n$, $P_{p,q}(y)$ is divisible by $Q(y)$. If $p$ is not divisible by $n$, then doing first the cancellation as above and then using the symmetry present in the sequence of remainders, the sum of the sequence of remainders equals 0 only if $r+1=t$, that is, only if $q$ is divisible by $n$.

2) We assume now that $n$ divides $p$ or $q$ and count the extra tiles $B$ that appears. They are counted by the coefficients of the quotient, call it $S(y)$, of the division of $P_{p,q}(y)$ by $Q(y)$. We need to compute the sum $S_1$ of the coefficients in $S(y)$ of the even powers of $y$  and the sum $S_2$ of the coefficients in $S(y)$ of the odd powers of $y$. The difference $S_1-S_2$ gives the number of extra tiles $B$ that we need to consider.

We use the equation relating the derivatives:
\begin{equation*}
\begin{gathered}
P'_{p,q}(y)=Q'(y)S(y)+Q(y)S'(y).
\end{gathered}
\end{equation*}
Note that $Q(-1)=0, Q'(-1)=n/2, S(-1)=S_1-S_2$. Plugging in $x=-1$  gives:
\begin{equation*}
\begin{gathered}
S_1-S_2=S(-1)=\frac{2P'_{p,q}(-1)}{n}.
\end{gathered}
\end{equation*}

Differentiating the equation of $P_{p,q}$ one has:
\begin{equation*}
\begin{gathered}
P'_{p,q}(y)=2\cdot 1+3\cdot 2y+4\cdot3y^2+\dots+(p-1)(p-2)y^{p-3}+p(p-1)y^{p-2}+\dots +p(q-1)y^{q-2}\\
+(p-1)qy^{q-1}+(p-2)(q+1)y^{q}+\dots +2(p+q-3)y^{p+q-4}+(p+q-2)y^{p+q-3}.
\end{gathered}
\end{equation*}

While computing $P_{p,q}(-1)$ we recall that $n$ is even and distinguish the following cases: \emph{Case A.} $p$ even, $q$ odd, \emph{Case B.} $p$ odd, $q$ even, \emph{Case C.} $p$ even, $q$ even.

We need the following formulas:
\begin{equation*}
\begin{gathered}
2\cdot 1-3\cdot 2+4\cdot3-\dots-(p-1)(p-2)=-\frac{p(p-2)}{2}\\
p(p-1)-p(p)+p(p+1)-\dots+p(q-2) -p(q-1)=-\frac{p(q-p+1)}{2}\\
(p-1)q-(p-2)(q+1)+\dots+3(p+q-4) -2(p+q-3)+(p+q-2)=\frac{pq}{2}.
\end{gathered}
\end{equation*}

\emph{Case A.} One has:
\begin{equation*}
\begin{gathered}
P'_{p,q}(-1)=2\cdot 1-3\cdot 2+4\cdot3-\dots-(p-1)(p-2)+p(p-1)-\dots -p(q-1)\\
+(p-1)q-(p-2)(q+1)+\dots -2(p+q-3)+(p+q-2)=\frac{p}{2}.
\end{gathered}
\end{equation*}

The number of extra $B$ tiles is $-\frac{p}{2}+\frac{p}{n}=\frac{p(1-k)}{n}.$ To have a complete reduction, the number of $B$ tiles has to be a multiple of $k-2$. As $k-1$ and $k-2$ are relatively prime, $p$ has to be a multiple of $n(k-2)$.

\emph{Case B.} One has:
\begin{equation*}
\begin{gathered}
P'_{p,q}(-1)=2\cdot 1-3\cdot 2+4\cdot3-\dots+(p-1)(p-2)-p(p-1)+\dots +p(q-1)\\
-(p-1)q+(p-2)(q+1)+\dots -2(p+q-3)+(p+q-2)=\frac{q}{2}.
\end{gathered}
\end{equation*}

The number of extra $B$ tiles is $\frac{q}{n}.$  We have the condition that $q$ is a multiple of $n(k-2)$.

\emph{Case C.} One has:
\begin{equation*}
\begin{gathered}
P'_{p,q}(-1)=2\cdot 1-3\cdot 2+4\cdot3-\dots-(p-1)(p-2)+p(p-1)-\dots +p(q-1)\\
-(p-1)q+(p-2)(q+1)+\dots +2(p+q-3)-(p+q-2)=0.
\end{gathered}
\end{equation*}

The number of extra $B$ tiles is $-\frac{p}{2}+\frac{p}{2}=0.$ In this case a signed tiling is always possible.

\section{Proof of Theorem~\ref{thm:main++}}\label{s:4-bis}

\begin{figure}[h!]
\centering
\begin{subfigure}{.15\textwidth}
\begin{tikzpicture}[scale=1]
\draw [line width = 1,<->] (0,2.5)--(0,0)--(2.5,0);
\draw [line width = 1] (0,2)--(2,2)--(2,0);
\draw [line width = 1] (1,0)--(1,2);
\draw [line width = 1] (0,1)--(2,1);

\node at (.5,.5) {$st$};
\node at (.5,1.5) {$st$};
\node at (1.5,.5) {$st$};
\node at (1.5, 1.5) {$st$};
\end{tikzpicture}
\caption{$p=2s, q=2t$}
\end{subfigure}
~~
\begin{subfigure}{.2\textwidth}
\begin{tikzpicture}[scale=1]
\draw [line width = 1,<->] (0,3.5)--(0,0)--(2.5,0);
\draw [line width = 1] (0,2)--(2,2)--(2,0);
\draw [line width = 1] (1,0)--(1,3);
\draw [line width = 1] (0,1)--(2,1);
\draw [line width = 1] (0,3)--(2,3)--(2,0);

\node at (.5,.5) {$s$};
\node at (.5,1.5) {$st$};
\node at (1.5,.5) {$s$};
\node at (1.5, 1.5) {$st$};
\node at (.5,2.5) {$st$};
\node at (1.5, 2.5) {$st$};
\end{tikzpicture}
\caption{$p=2s, q=2t+1$}
\end{subfigure}
~~
\begin{subfigure}{.2\textwidth}
\begin{tikzpicture}[scale=1]
\draw [line width = 1,<->] (0,2.5)--(0,0)--(3.5,0);
\draw [line width = 1] (0,2)--(2,2)--(2,0);
\draw [line width = 1] (1,0)--(1,2);
\draw [line width = 1] (0,1)--(3,1);
\draw [line width = 1] (2,2)--(3,2)--(3,0);

\node at (.5,.5) {$t$};
\node at (.5,1.5) {$t$};
\node at (1.5,.5) {$st$};
\node at (1.5, 1.5) {$st$};
\node at (2.5,.5) {$st$};
\node at (2.5, 1.5) {$st$};
\end{tikzpicture}
\caption{$p=2s+1, q=2t$}
\end{subfigure}
~~~~
\begin{subfigure}{.2\textwidth}
\begin{tikzpicture}[scale=1]
\draw [line width = 1,<->] (0,3.5)--(0,0)--(3.5,0);
\draw [line width = 1] (0,2)--(2,2)--(2,0);
\draw [line width = 1] (1,0)--(1,3);
\draw [line width = 1] (0,1)--(3,1);
\draw [line width = 1] (2,2)--(3,2)--(3,0);
\draw [line width = 1] (0,3)--(3,3)--(3,0);
\draw [line width = 1] (2,2)--(2,3);

\node at (.5,.5) {$1$};
\node at (.5,1.5) {$t$};
\node at (1.5,.5) {$s$};
\node at (1.5, 1.5) {$st$};
\node at (2.5,.5) {$s$};
\node at (2.5, 1.5) {$st$};
\node at (.5, 2.5) {$t$};
\node at (1.5,2.5) {$st$};
\node at (2.5, 2.5) {$st$};
\end{tikzpicture}
\caption{$p=2s+1, q=2t+1$}
\end{subfigure}
\caption{Reduced configurations.}
\label{fig:newcong83-new01}
\end{figure}
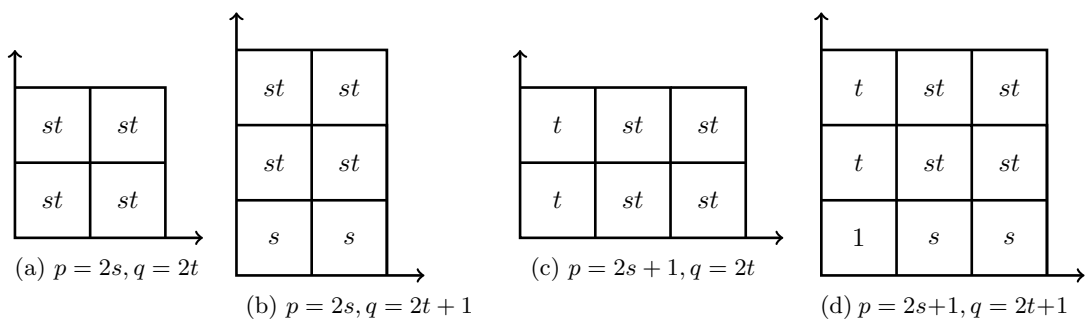

Let $k\ge 3$. Consider a $q\times p, q\ge p\ge 1,$ rectangle. Using the presence of $D_1(k)$ and $\bar D_1(k)$ in the ideal, the rectangle can be reduced to one of the configurations in Figure~\ref{fig:newcong83-new01}. The configuration in a), and its copies appearing in b), c), d), are multiples of $H_5(k)$ and can be reduced to zero. The remaining region in b) can be reduced to $s[1-(k-2)(k-1)](1+y)$. As $1-(k-2)(k-1)$ is never a multiple of $k(k-2)$ for $k\ge 4$, this configuration can be reduced further by $D_2(k)$ to zero only if $s$ is a multiple of $k(k-2)$. Same reasoning works for c). The remaining region in d) can be reduced further to $[(k-2)(k-1)s-t]y+[(k-2)(k-1)s+s-t+1]$ which is never a multiple of $D_2(k)$, thus cannot be reduced to zero.

Assume now $k=2$. The proof is similar, but one observes that only configuration a) in Figure~\ref{fig:newcong83-new01} can be reduced to zero using the Gr\"obner basis for $\mathcal{T}_4^+$.

\section{Proof of Proposition~\ref{thm:main-coro93}}\label{s:5}

If $k$ is even, finding a signed tiling for a $k$-inflated copy of the $L$ $n$-omino can be reduced, via reductions by $C_3(k), C_4(k)$ tiles, to finding a signed tiling for a $nk\times k$ rectangle. From Theorem~\ref{thm:main} follows that such a tiling always exists.
If $k$ is odd, a reduction to a $kn\times k$ rectangle can be done only modulo a $B$ tile, which does not belong to the ideal generated by $\mathcal{T}_n$.

\section{The method of Barnes for $\mathcal{T}_n, n$ even}\label{s:barnes}\label{s:6}

In this section we give a proof of Theorem~\ref{thm-main-barnes} following a method of Barnes. We assume familiarity with \cite{Barnes1, Barnes2}. We apply the method to $\mathcal{T}_n, n\ge 6$ even. Consider the polynomials  \eqref{eq:generators-k-even} associated to the tiles in $\mathcal{T}_n$ and denote by $I$ the ideal generated by them. We show that the algebraic variety $V\subset \mathbb{C}^2$ defined by  \eqref{eq:generators-k-even} consists only of the points $(\epsilon, 1+\epsilon+\epsilon^2+\dots \epsilon^{n-1}),$ where $\epsilon$ is an $n$-th root of identity different from 1.

Separating $x$ from $H_1(k)=0$, replacing in $H_2(k)=0$ and factoring the resulting polynomial gives:
\begin{equation*}
(y^{2k-1}+y^{2k-2}+\dots +y^2+y+1)(y^{2k-3}+y^{2k-4}+\dots +y^2+y+1)=0.
\end{equation*}

Denote the polynomial on the left hand side by $f_2(y)$, and denote the corresponding polynomial in the variable $x$ (obtained from $H_3(k)$ and $H_4(k)$) by $f_1(x)$. Their roots are roots of unity of order $2k-1$ and $2k-3$. Using the system of equations that defines $V$, the roots of order $2k-3$ can be eliminated. Moreover, the only solutions of the system are as above.

We show now that $I$ is a radical ideal. We use an algorithm of Seidenberg which can be applied to find the radical ideal of a zero dimensional algebraic variety over an algebraically closed field. See Lemma 92 in~\cite{Seidenberg}. Compare also with \cite[Theorem 7.1]{Barnes1}. As $V\subset \mathbb{C}^2$ is zero dimensional, one can find square free polynomials $\bar f_1(x)$ and $\bar f_2(y)$ that belong to the radical ideal. We take these to be $f_1(x)=F(x), f_2(y)=F(y),$ where:
\begin{equation*}
\begin{gathered}
F(x)=(k-2)(x^{2k-1}+x^{2k-2}+\dots +x^2+x+1).
\end{gathered}
\end{equation*}
The ideal generated by $H_i(k)'$s and $\bar f_1(x), \bar f_2(y)$ is radical. To show that $I$ is radical, it is enough to show that $\bar f_1(x), \bar f_2(y)$ belong to $I$. For $\bar f_1(x)$ we have $\bar f_1(x)=(k-2)xH_3(k)+C_5(k)$ and $\bar f_2(y)$ follows by symmetry.

We apply now the main result in \cite[Lemma 3.8]{Barnes1}: a region $R$ is signed tiled by $\mathcal{T}_n$ if and only if the polynomial $f_R(x,y)$ associated to $R$ evaluates to zero in all point of $V$. If $R$ is a $p\times q$ rectangle, then $f_R(x,y)=\frac{x^q-1}{x-1}\cdot \frac{y^p-1}{y-1},$
which evaluates to zero in all points of $V$ if and only if one of $p, q$ is divisible by $n$. This gives the first statement in Theorem~\ref{thm-main-barnes}.

For the second statement in Theorem~\ref{thm-main-barnes} we use the method described in the proof of \cite[Theorem 4.2]{Barnes1}. A set of generators over $\mathbb{Q}$ for the rectangles that have a side divisible by $n$ is given by the set $\mathcal H=\{H_1(k), H_2(k),H_3(k), H_4(k)\}$ and the polynomial with rational coefficients $\frac{1}{k-2}C_5(k)$. As $C_5(k)$ is already generated by $\mathcal H$, this implies that $(k-2)$ multiples of the elements in $\mathcal{H}$ can signed tile with integer coefficients any $(k-2)$ multiple of a rectangle with a side divisible by $n$.

\section{The method of Barnes for $\mathcal{T}^+_n, n$ even}\label{s:barnes+}

Let $k\ge 3$. Adding the extra polynomial $1+x+y+xy$ to the set of generators, reduces the variety $V$ to the point $(-1,-1)$. As before, the ideal $I$ is radical and the square free polynomials $\bar f_1, \bar f_2$ can be chosen to be $\bar f_1(x)=k(k-2)(x+1), \bar f_2(x)=k(k-2)(y+1)$. The second one belongs to the Gr\"obner basis for $\mathcal{T}^+_n$ and the first one can be generated as well as our set of generators $H_i(k)$ is symmetric in the variables $x,y$. The statement in Theorem~\ref{t:barnes+} follows now from~\cite[Lemma 3.8]{Barnes1}.

Assume now $k=2$. In this case the ideal $I$ it is not radical. This follows using the theory developed in~\cite{Barnes1} about colorings. It is shown in~\cite{Barnes1} that $\mathcal{T}_4$ has 4 colorings, three standard and one nonstandard due to the differential operator $\partial_x+\partial_y$. It is easy to check that one of the standard colorings~\cite[Figure 1]{Barnes1} and the nonstandard coloring~\cite[Figure 4]{Barnes1} are the only colorings for $\mathcal{T}_4^+$. One can check that a rectangle fits these colorings only if and only if it has both sides even, so it follows from~\cite[Theorem 5.3]{Barnes1} that a rectangle is signed tiled by $\mathcal{T}_4^+$ if and only if it has both sides even.



\section*{Acknowledgement}

V. Nitica was partially supported by Simons Foundation Grant 208729. While working on this project, K. Gill was undergraduate student at West Chester University.

\end{document}